\documentclass[11pt]{article}
\usepackage[a4paper,margin=1in]{geometry}
\usepackage{amsmath,amssymb,amsthm,mathtools}
\usepackage{extarrows}
\usepackage{graphicx,subcaption,booktabs,array}
\graphicspath{{./}{Figure/}{../Figure/Figures_Aug_20th_2026/}{Figure/Figures_Aug_20th_2026/}}
\usepackage{microtype}
\usepackage{comment}
\usepackage{enumitem}

\usepackage{xcolor}
\usepackage{tikz}
\usetikzlibrary{calc}
\definecolor{gtblue}{RGB}{0,114,178}
\definecolor{gtorange}{RGB}{213,94,0}
\definecolor{gtgreen}{RGB}{0,158,115}
\definecolor{gtlabel}{RGB}{70,90,140}

\usepackage[hidelinks]{hyperref}
\newtheorem{theorem}{Theorem}[section]
\newtheorem{proposition}[theorem]{Proposition}
\newtheorem{lemma}[theorem]{Lemma}

\theoremstyle{definition}
\newtheorem{definition}[theorem]{Definition}
\newtheorem*{definition*}{Definition}
\theoremstyle{remark}
\newtheorem{remark}[theorem]{Remark}
\newcommand{\RT}{\mathbb R^\Gamma_T}
\newcommand{\supp}{\operatorname{supp}}
\newcommand{\smin}{\sigma_{\min,T}}
\DeclareMathOperator{\diag}{diag}
\DeclareMathOperator{\Var}{Var}
\title{Stability of Sparse Crystallographic Phase Retrieval}
\author{Ruilin Zhang \and Hai Zhang}
\date{}
\begin{document}
\maketitle
\begin{abstract}
Crystallographic phase retrieval seeks to recover a periodic signal from its
Fourier intensities. The problem is known to be ill-posed without additional constraints. This paper investigates how sparsity can restore stability.
For real signals on a finite periodic grid of size \(N\), we prove local
Lipschitz stability at generic signals on a random support at
the \(N^{1-\varepsilon}\) scale, with high probability. We further establish quantitative local
Lipschitz stability with high probability for sparsity up to order \(\sqrt N\).
We also reveal that nonnegativity can upgrade local stability to global stability,
which may help explain its effectiveness as a prior in phase retrieval.
Specifically, under the same conditions as our quantitative local stability result,
adding a nonnegativity constraint yields global Lipschitz stability for
signals with sparsity order $\sqrt{N}$.
Our quantitative results rely on autocorrelation entries to which exactly one
pair of support points contributes. These pairs form the edges of a graph $G$
on the support, called the graph of isolated autocorrelations, and a lower bound on the smallest eigenvalue of the signless
Laplacian of $G$ yields an upper bound on the Lipschitz stability constant.
Numerical experiments suggest that this graph-based bound on the Lipschitz constant degrades beyond the
\(\sqrt N\) sparsity scale, while local stability may persist up to sparsities
near \(N/2\).
\end{abstract}

\section{Introduction}
\label{sec:introduction}

Phase retrieval seeks to reconstruct a signal from the magnitudes of linear
measurements. Its classical Fourier form arises when a detector records
diffraction intensities but not Fourier phases. In X-ray crystallography,
these intensities carry information about a periodic electron density
\cite{Millane1990}. Related problems occur in astronomical imaging, electron
microscopy, and coherent diffraction imaging
\cite{DaintyFienup1987,Zuo2003,Miao1999}. They have motivated reconstruction methods for a variety of sampling models and
signal constraints where phase information is not available.
Phase retrieval has also attracted considerable attention in applied
mathematics \cite{Shechtman2015Overview,Jaganathan2015Overview}. As an inverse problem, it is simple to state yet nonlinear, nonconvex \cite{Luke2017WhatsNew}, and ill-posed \cite{BendoryEdidin2020}.

In a general $d$-dimensional formulation with \(d\geq1\), the signals are indexed by
\(\Gamma=\mathbb Z_{n_1}\times\cdots\times\mathbb Z_{n_d}\). Write \(N=|\Gamma|\).
Let \(x\in\mathbb C^\Gamma\) (or \(x\in\mathbb R^\Gamma\) in the real case)
be the unknown signal. Given known measurement vectors
\(a_0,\ldots,a_{m-1}\in\mathbb C^\Gamma\), we observe the intensity measurements
\begin{equation}
 y_k=|\langle a_k,x\rangle|^2,\qquad k=0,\ldots,m-1.
 \label{eq:measurement-model}
\end{equation}
Common choices of measurement vectors include independent real or complex
Gaussian random vectors~\cite{Candes2013} and Fourier vectors.

In crystallographic phase retrieval, signals are periodic, and \(\Gamma\) is
endowed with its natural abelian group structure. The measurement vectors are
the Fourier vectors on \(\Gamma\). Equation~\eqref{eq:measurement-model} then takes the form
\[
  y = |\mathcal{F}(x)|^2,
\]
where \(\mathcal{F}\) denotes the unmasked unitary discrete Fourier transform
on \(\Gamma\) and the squared modulus is taken entrywise. Throughout this paper,
we consider real signals \(x \in \mathbb{R}^\Gamma\).

For real \(x\), \(|\mathcal F(x)|^2\) is invariant under global sign change,
cyclic translation \(x_j\mapsto x_{j-v}\) with \(v\in\Gamma\), and reflection
\(x_j\mapsto x_{-j}\), where all indices are taken in \(\Gamma\); these are the
trivial ambiguities~\cite{BendoryEdidin2020}. Unless otherwise stated,
uniqueness and recovery in the Fourier setting are understood modulo these
ambiguities.

Algorithms for Fourier phase retrieval have been developed for several measurement and constraint models. The Gerchberg--Saxton algorithm, introduced in 1972, alternates between
the image and diffraction planes, enforcing the measured modulus in
each plane \cite{Gerchberg1972}. Fienup's error-reduction and input--output algorithms work with a single Fourier-magnitude measurement and combine it with object-domain constraints such as support and nonnegativity. The hybrid input--output (HIO) variant adds a feedback term to the object-domain update \cite{Fienup1982}. Later iterative projection algorithms include the difference map (DM) \cite{elser2003phase} and relaxed averaged alternating reflections (RAAR) \cite{Luke2005}. Marchesini compares these methods in a unified framework \cite{Marchesini2007}.

For generalized phase retrieval with random measurements, such as independent Gaussian measurement vectors, convex relaxations like PhaseLift and PhaseCut \cite{Candes2013,Waldspurger2015} recover the signal exactly up to a global phase with high probability from sufficiently many measurements. However, these guarantees depend on the randomness of the measurement vectors and do not carry over directly to the Fourier phase retrieval problem we consider here.

\subsection{Sparsity and well-posedness}
\label{subsec:sparsity}

Phase retrieval algorithms typically use prior information to make the
reconstruction problem well-posed. Sparsity is a common prior: it requires \(|\supp(x)|\leq s\) in a
fixed basis. 
Shechtman, Beck, and Eldar's GESPAR combines support swaps with damped 
Gauss--Newton optimization on the current support
\cite{Shechtman2014GESPAR}. Some of their numerical experiments suggest
reliable reconstruction at sparsity scale \(N^{1/3}\).
For noiseless oversampled aperiodic Fourier data, the two-stage TSPR algorithm recovers random sparse signals with probability tending to one in the regime
\(s=O(N^{1/2-\varepsilon})\), for fixed
\(0<\varepsilon<\frac12\)
\cite[Theorem~III.1]{jaganathan2017sparse}.
For periodic crystallographic benchmarks, Elser, Lan, and Bendory
\cite{ElserLanBendory2018Benchmark} observe that the mean iteration
count of relaxed-reflect-reflect (RRR) grows approximately exponentially
with \(s^2/V\), where \(s\) counts atoms and \(V\) is the effective
number of Fourier samples. This suggests keeping \(s=O(\sqrt V)\) for practical computation (\(\sqrt N\) when
\(V\asymp N\)).

The empirical success of these algorithms raises a natural question: how does
sparsity of the signal affect the well-posedness of Fourier phase retrieval,
that is, the uniqueness and stability of the reconstruction? This question
is the main focus of this paper. We first review several representative results
in this direction.

For one-dimensional aperiodic Fourier data, Ranieri et al.
\cite[Theorem~2]{Ranieri2013} prove uniqueness for real collision-free
\(s\)-sparse signals with \(s\geq3\) and \(s\neq6\), where
collision-free means that all pairwise differences of support points are
distinct.
For periodic Fourier data, Edidin and Suresh
\cite[Theorem~1.1(ii)]{EdidinSuresh2026} prove that a generic real
\(s\)-sparse signal in a generic basis is uniquely determined among
\(s\)-sparse real signals up to sign when
\(s\leq\frac N2\).
In the standard coordinate basis on \(\mathbb Z_N\), Bendory and
Edidin \cite[Conjecture~1]{BendoryEdidin2020} conjecture generic
uniqueness among real \(s\)-sparse signals when \(s\leq\frac N2\).

Uniqueness alone does not ensure reliable reconstruction; stability is also
important in practice: numerical experiments~\cite{barnett2020geometry} suggest a link
between poor local stability and stagnation, the phenomenon in which further
iterations make little or no progress toward a
solution~\cite{fienup1986stagnation,elser2003phase}. For suitable frame
measurements, global stability has been
established~\cite{Bandeira2014SavingPhase,BalanZou2016,EldarMendelson2014}. For
fixed, unmasked Fourier measurements, however, stability is usually studied
locally and there are few results.

For real periodic signals, the analysis of Ghosh and Rigollet
\cite{ghosh2023sparse} can be used to derive local Lipschitz stability
for collision-free signals.
The local estimate compares nearby signals on the same support.
Their analysis also yields local Lipschitz stability within a symmetric
sparse class, where \(x_j=x_{-j}\). In particular, this holds with high probability
for signals sampled from the symmetric Bernoulli--Gaussian model with
sparsity parameter \(\log^9 N\leq s\leq N/\log^5 N\).

\subsection{Related work}
\label{subsec:related-work}

In addition to the local stability results for sparse signals
\cite{ghosh2023sparse} discussed in the preceding subsection, we also
refer to work on local stability in the broader setting of Fourier
phase retrieval. Barnett, Epstein, Greengard, and Magland
\cite{barnett2020geometry,barnett2022geometry} relate the
well-posedness of reconstruction to the geometry of the measurement and
support constraints.

More specifically, for a real image \(f\), write
\(a=\mathcal M(f)=|\widehat f|\), where \(\widehat f\) is its
discrete Fourier transform. Denote by
\(A_a=\{g:\mathcal M(g)=a\}\) the magnitude torus and by \(B\)
the linear space of real images satisfying the support constraint.  Assume that all Fourier coefficients of $f$
are nonzero, so that $A_a$ is a smooth torus. They
prove that the \emph{transversality} condition
\[
 T_f A_a\cap B=\{f\},
\]
where \(T_f A_a\) is the affine tangent space to \(A_a\) through
\(f\), is equivalent to the existence of constants
\(C,\eta>0\) such that
\[
 C\|\delta f\|_2
 \leq\|\mathcal M(f+\delta f)-\mathcal M(f)\|_2,
 \qquad\delta f\in B,\quad\|\delta f\|_2<\eta.
\]
They also identify some cases in which this local Lipschitz condition fails.

\subsection{Contributions and main results}
\label{subsec:main-results}

For a nonempty sparse support constraint \(T\subseteq\Gamma\), define
\[
 \RT=\{x\in\mathbb R^\Gamma:\supp(x)\subseteq T\}.
\]
We study the phase retrieval problem of reconstructing a ground-truth signal
\(x^*\in\RT\) from its Fourier intensity \(|\mathcal F(x^*)|^2\) among signals
in \(\RT\).

We consider two notions of stability.

\begin{itemize}
\item \emph{Local stability.}
We say that stability holds locally at \(x^*\) if there exist a constant
\(L_{\mathrm{local}} = L_{\mathrm{local}}(x^*) \in (0,\infty)\) and a
neighborhood \(\mathcal U \subseteq \RT\) of \(x^*\) such that
\[
  \|x^{(1)} - x^{(2)}\|_2
  \leq L_{\mathrm{local}}
  \bigl\| \, |\mathcal F(x^{(1)})|^2 - |\mathcal F(x^{(2)})|^2 \, \bigr\|_2
  \qquad \text{for all } x^{(1)}, x^{(2)} \in \mathcal U.
\]

\item \emph{Global stability.}
We say that stability holds globally on a signal class \(\mathcal K \subseteq \RT\)
if the above inequality holds for all \(x^{(1)}, x^{(2)} \in \mathcal K\)
with a constant
\(L_{\mathrm{global}} = L_{\mathrm{global}}(\mathcal K) \in (0,\infty)\)
that depends only on \(\mathcal K\).
\end{itemize}
We say that \emph{generic local stability} holds on \(\RT\) if local stability
holds at Lebesgue-almost every \(x^*\in\RT\).

Our main results are as follows.

\paragraph{Generic local stability at the \(N^{1-\varepsilon}\) scale.}
Fix \(\varepsilon\in(0,1]\) and \(C>0\), and let \(t=t(N)\) be
an integer satisfying\footnote{For \(u\in\mathbb R\), \(\lceil u\rceil\)
denotes the smallest integer greater than or equal to \(u\).}
\[
 (2+2^d)\lceil\varepsilon^{-1}\rceil\leq t\leq CN^{1-\varepsilon}.
\]
Let \(T\) be a uniformly random \(t\)-element subset of \(\Gamma\).
Theorem~\ref{thm:generic} shows that, for all sufficiently large \(N\),
with probability at least \(1-O(N^{-\varepsilon})\), generic local
stability holds on \(\RT\).

This result shows that instability is exceptional under random sparse
support constraints. It complements the work of Barnett et
al.~\cite{barnett2020geometry}, who identify specific configurations in
which local Lipschitz stability fails.

\paragraph{Quantitative local stability at the \(\sqrt{N}\) scale.}
Fix \(0<\alpha\leq1\) and \(0<m\leq M\). Let \(T\subseteq\Gamma\) be a
Bernoulli random support with \(\mathbb E|T|=t\). That is, each
point of \(\Gamma\) is included in \(T\) independently with probability \(t/N\),
where
\begin{equation}
  (\log N)^4=o(t),\qquad t\leq c(\alpha)\sqrt N
  \label{eq:t-range}
\end{equation}
for a constant \(c(\alpha)>0\) depending only on \(\alpha\).
Theorem~\ref{thm:random-stability} shows that, for all sufficiently large
\(N\), with probability at least
$
  1-\exp\!\bigl[-\Omega\bigl(\sqrt t\,\log(N/t)\bigr)\bigr],
$\footnote{We write \(f=\Omega(g)\) if there is a constant \(c>0\) such
that \(f\geq c\,g\) for all sufficiently large \(N\).}
the random support \(T\) guarantees
\[
  L_{\mathrm{local}}(x^*)\leq O\!\left(\frac{M}{m^2}\sqrt{\frac{N}{\alpha t}}\right)
\]
for every signal
\[
  x^*\in\Bigl\{x\in\RT:\ |\supp(x)|\geq\alpha|T|,\ \
  m\leq|x_i|\leq M\ \text{for all } i\in\supp(x)\Bigr\}.
\]

This extends the collision-free local stability estimates of Ghosh and
Rigollet \cite{ghosh2023sparse}. The collision-free condition holds with high probability
for random supports of size \(o(N^{1/4})\).
Moreover, their local estimates compare signals on the same support,
corresponding to the exact-support case \(\supp(x^*)=T\) in our formulation.

We obtain the bound on \(L_{\mathrm{local}}(x^*)\) from the pairs
\(\{i,j\}\subseteq T\) such that no other pair in \(T\) has difference
\(\pm(i-j)\). For such a pair, the autocorrelation of \(x\in\RT\) at
\(i-j\) depends only on \(x_ix_j\). These pairs form the edges of a graph on
\(T\), the graph of isolated autocorrelations
(Subsection~\ref{subsec:isolated-graph}). The bound then follows from a
lower bound on the smallest eigenvalue of its signless Laplacian.
Jaganathan, Oymak, and Hassibi~\cite{jaganathan2017sparse} used a related
graph for aperiodic phase retrieval.

\paragraph{Global stability at the \(\sqrt N\) scale for nonnegative signals.}
Let \(T\subseteq\Gamma\) be a Bernoulli random support with
\(\mathbb E|T|=t\), where \(t\) satisfies \eqref{eq:t-range}.
Theorem~\ref{thm:random-positive-global} shows that, for all sufficiently
large \(N\), with probability at least
$
  1-\exp\!\bigl[-\Omega\bigl(\sqrt t\,\log(N/t)\bigr)\bigr],
$
the support \(T\) admits global stability with
\[
 L_{\mathrm{global}}
 \leq O\!\left(\frac{M}{m^2}\sqrt{\frac{N}{\alpha t}}\right)
\]
on
\[
 \mathcal K=\left\{x\in\RT\cap\mathbb R^\Gamma_+:
 |\supp(x)|\geq\alpha|T|,\quad
 m\leq x_i\leq M\ \text{for }i\in\supp(x)\right\},
\]
where \(\mathbb R^\Gamma_+\) is the nonnegative orthant.

\section{Autocorrelation and  stability}
\label{sec:autocorrelation}

We introduce the autocorrelation and use it to investigate the stability of phase retrieval. In contrast to the local tangent-space method of Barnett et al.~\cite{barnett2020geometry}, the quadratic structure of the autocorrelation yields
an exact midpoint identity that links local and global stability estimates.

\subsection{Autocorrelation and its Jacobian}

\begin{definition}[Circular autocorrelation]
\label{def:autocorrelation}
The circular autocorrelation map \(A:\mathbb R^\Gamma\to\mathbb R^\Gamma\)
is defined by
\[
 A(x)_v=\sum_{i\in\Gamma}x_ix_{i-v},
 \qquad x\in\mathbb R^\Gamma,\quad v\in\Gamma,
\]
where, for \(i=(i_1,\dots,i_d)\) and \(v=(v_1,\dots,v_d)\) in \(\Gamma\),
\[
 i\pm v=\bigl((i_1\pm v_1)\bmod n_1,\ \dots,\ (i_d\pm v_d)\bmod n_d\bigr).
\]
All sums and differences of indices in this paper are computed in this way.
\end{definition}

Note that \(\mathcal F\) is unitary, i.e.,
\(\|\mathcal F(x)\|_2=\|x\|_2\), and that
\(\mathcal F(A(x))=\sqrt N\,|\mathcal F(x)|^2\). It follows
\begin{equation}
 \|A(x^{(1)})-A(x^{(2)})\|_2
 =\sqrt N\,\bigl\|\,|\mathcal F(x^{(1)})|^2-|\mathcal F(x^{(2)})|^2\,\bigr\|_2.
 \label{eq:normalization}
\end{equation}

\medskip
\begin{definition}[Jacobian of \(A\)]
\label{def:jacobian}
For \(x\in\mathbb R^\Gamma\), let
\(A'(x)\in\mathbb R^{\Gamma\times\Gamma}\) denote the Jacobian of
\(A\) at \(x\). After ordering the elements of \(\Gamma\), it can be viewed
as an \(N\times N\) real matrix with entries
\begin{equation}
 A'(x)_{v,i}=\frac{\partial A(x)_v}{\partial x_i}
 =x_{i-v}+x_{i+v}=A'(x)_{-v,i},\qquad v,i\in\Gamma.
 \label{eq:jacobian}
\end{equation}
\end{definition}

For \(x,h\in\RT\) and \(v\in\Gamma\), \eqref{eq:jacobian} gives
\begin{equation}
 (A'(x)h)_v
 =\sum_{i\in\Gamma}x_{i-v}h_i+\sum_{j\in\Gamma}x_{j+v}h_j
 =\sum_{\substack{i,j\in T\\i-j=v}}x_jh_i
  +\sum_{\substack{i,j\in T\\i-j=v}}x_ih_j
 =\sum_{\substack{i,j\in T\\i-j=v}}(x_jh_i+x_ih_j).
 \label{eq:Ah-pairs}
\end{equation}

Denote
\[
\Gamma[2]=\{v\in\Gamma:v=-v\}=\{v\in\Gamma:2v=0\}.
\]
An element \(v\in\Gamma=\mathbb Z_{n_1}\times\cdots\times\mathbb Z_{n_d}\)
satisfies \(2v=0\) if and only if each of its coordinates does. In
\(\mathbb Z_n\), the solutions of \(2u=0\) are \(u=0\) and, for even \(n\),
\(u=n/2\). Hence \(|\Gamma[2]|\leq2^d\).

By \eqref{eq:jacobian}, the rows of \(A'(x)\) indexed by \(v\) and \(-v\)
coincide, so \(A'(x)\) has at most one distinct row for each set
\(\{v,-v\}\), \(v\in\Gamma\). Hence
\begin{equation}
 \operatorname{rank}A'(x)
 \leq|\Gamma[2]|+\frac{N-|\Gamma[2]|}{2}
 =\frac{N+|\Gamma[2]|}{2}.
 \label{eq:row-bound}
\end{equation}

With \(A\) and \(A'\) defined above, we note that the quadratic structure of
\(A\) gives the following identity.

\begin{lemma}[Midpoint identity]
\label{lem:midpoint}
For any \(x^{(1)},x^{(2)}\in\mathbb R^\Gamma\),
\begin{equation}
 A(x^{(1)})-A(x^{(2)})=A'\left(\frac{x^{(1)}+x^{(2)}}{2}\right)(x^{(1)}-x^{(2)}).
 \label{eq:midpoint}
\end{equation}
\end{lemma}

\subsection{Jacobian criteria for local and global stability}
\label{subsec:stability-criteria}

We first characterize local Lipschitz stability at a signal \(x^*\), and
then give a sufficient condition for global Lipschitz stability on a
signal class \(\mathcal K\). We state both bounds in terms of \(A(x)\); by
\eqref{eq:normalization}, they can be rewritten in terms of
\(|\mathcal F(x)|^2\). 

Recall that
\(\RT=\{x\in\mathbb R^\Gamma:\supp(x)\subseteq T\}\). For a matrix \(\Phi\)
with columns indexed by \(\Gamma\) and a symmetric matrix
\(H\in\mathbb R^{\Gamma\times\Gamma}\), define
\[
 \sigma_{\min,T}(\Phi)=\min_{\substack{h\in\RT\\\|h\|_2=1}}\|\Phi h\|_2,
 \qquad
 \lambda_{\min,T}(H)=\inf_{\substack{h\in\RT\\\|h\|_2=1}}h^{\mathsf T}Hh.
\]

\par\pagebreak[0]
\begin{samepage}
\begin{lemma}[Local Lipschitz stability]
\label{lem:local}
Let \(\varnothing\neq T\subseteq\Gamma\) and \(x^*\in\RT\). The following
are equivalent:
\begin{enumerate}[label=(\roman*)]
\item \(A'(x^*)\) is injective on \(\RT\).
\item There exist a neighborhood \(\mathcal U\) of \(x^*\) in \(\RT\) and
\(0<L<\infty\) such that, for all \(x^{(1)},x^{(2)}\in\mathcal U\),
\[
 \|x^{(1)}-x^{(2)}\|_2\leq L\|A(x^{(1)})-A(x^{(2)})\|_2.
\]
\end{enumerate}
If \(\sigma=\smin(A'(x^*))>0\), one may take \(L=\frac{2}{\sigma}\)
and a convex neighborhood.
\end{lemma}
\end{samepage}
\begin{proof}
We first show that (i) implies (ii). Injectivity of \(A'(x^*)\) on \(\RT\)
is equivalent to \(\sigma=\smin(A'(x^*))>0\). Since \(A'\) depends
continuously on \(x\), we can choose a convex neighborhood \(\mathcal U\)
of \(x^*\) in \(\RT\) such that
\[
 \smin(A'(x))
 \geq\sigma-\|A'(x)-A'(x^*)\|
 \geq\frac\sigma2,
 \qquad x\in\mathcal U.
\]
For \(x^{(1)},x^{(2)}\in\mathcal U\), \((x^{(1)}+x^{(2)})/2\in
\mathcal U\) by convexity, so \eqref{eq:midpoint} gives
\[
 \|A(x^{(1)})-A(x^{(2)})\|_2\geq\frac\sigma2\|x^{(1)}-x^{(2)}\|_2.
\]
This proves (ii) with \(L=\frac{2}{\sigma}\).

We next show that (ii) implies (i). For \(h\in\RT\), applying the
estimate in (ii) to \(x^{(1)}=x^*+uh\) and \(x^{(2)}=x^*\), we obtain
\[
 \|h\|_2
 \leq L\lim_{u\to0}\frac{\|A(x^*+uh)-A(x^*)\|_2}{|u|}
 =L\|A'(x^*)h\|_2,
\]
so \(A'(x^*)h=0\) forces \(h=0\), and \(A'(x^*)\) is injective on \(\RT\).
\end{proof}

\begin{remark}
\label{rem:support-rank}
If \(|T|>(N+|\Gamma[2]|)/2\), then \eqref{eq:row-bound} gives
\(\operatorname{rank}A'(x)<|T|=\dim\RT\), so \(A'(x)\) is not injective
on \(\RT\), and hence
\[
 \smin(A'(x))=0,\qquad x\in\mathbb R^\Gamma.
\]
\end{remark}

\begin{remark}
\label{rem:fourier-geometry}
Set \(\mathcal M(x)=|\mathcal F(x)|\). Note that
\[
 \frac1{\sqrt N}\mathcal F(A'(x))
 =2\diag(|\mathcal F(x)|)\,\mathcal M'(x)
\]
whenever all Fourier coefficients of \(x\) are nonzero. In this case,
\(A'(x)\) and \(\mathcal M'(x)\) have the same kernel, so the results
in Lemma~\ref{lem:local} are equivalent to those
of \cite{barnett2020geometry} with magnitude measurements.
Our results for intensity measurements, however, do not require all
Fourier coefficients to be nonzero.
\end{remark}

\begin{lemma}[Global injectivity and stability]
\label{lem:global}
Let \(\varnothing\neq T\subseteq\Gamma\) and \(\mathcal K\subseteq\RT\).
Suppose that \(A'((x^{(1)}+x^{(2)})/2)\) is injective on \(\RT\) for
every pair \(x^{(1)},x^{(2)}\in\mathcal K\). Then \(A|_{\mathcal K}\)
is injective, and for all \(x^{(1)},x^{(2)}\in\mathcal K\),
\begin{equation}
 \|x^{(1)}-x^{(2)}\|_2
 \leq\frac{1}{\smin\Bigl(A'\Bigl(\frac{x^{(1)}+x^{(2)}}{2}\Bigr)\Bigr)}\,
       \|A(x^{(1)})-A(x^{(2)})\|_2.
 \label{eq:global-pairwise}
\end{equation}
If there exists \(\sigma>0\) such that
\(\smin(A'((x^{(1)}+x^{(2)})/2))\geq\sigma\) for all such pairs, then, for
all \(x^{(1)},x^{(2)}\in\mathcal K\),
\begin{equation}
 \|x^{(1)}-x^{(2)}\|_2
 \leq\frac1\sigma\|A(x^{(1)})-A(x^{(2)})\|_2.
 \label{eq:global-single-sigma}
\end{equation}
\end{lemma}
\begin{proof}
We fix \(x^{(1)},x^{(2)}\in\mathcal K\). Since \(\mathcal K\subseteq\RT\),
we have \(x^{(1)}-x^{(2)}\in\RT\), so \eqref{eq:midpoint} gives
\[
 \|A(x^{(1)})-A(x^{(2)})\|_2
 =\Bigl\|A'\Bigl(\frac{x^{(1)}+x^{(2)}}{2}\Bigr)(x^{(1)}-x^{(2)})\Bigr\|_2
 \geq\smin\Bigl(A'\Bigl(\frac{x^{(1)}+x^{(2)}}{2}\Bigr)\Bigr)
      \|x^{(1)}-x^{(2)}\|_2.
\]
By assumption, \(A'((x^{(1)}+x^{(2)})/2)\) is injective on \(\RT\), so
\(\smin(A'((x^{(1)}+x^{(2)})/2))>0\), and we obtain
\eqref{eq:global-pairwise}. In particular, \(A(x^{(1)})=A(x^{(2)})\) forces
\(x^{(1)}=x^{(2)}\), so \(A|_{\mathcal K}\) is injective.

If \(\smin(A'((x^{(1)}+x^{(2)})/2))\geq\sigma\) for all pairs in
\(\mathcal K\), then \eqref{eq:global-pairwise} yields
\eqref{eq:global-single-sigma}.
\end{proof}

\section{Generic local stability for random supports}
\label{sec:generic}

In this section, \(T\) is a uniformly random \(t\)-element subset of
\(\Gamma\): each \(t\)-element subset of \(\Gamma\) is sampled with
probability \(\binom Nt^{-1}\). We show that, for
\(t\leq CN^{1-\varepsilon}\), with high probability, local Lipschitz
stability holds at almost every signal in \(\RT\).

\begin{definition*}[Nonsingular support]
A nonempty support \(T\subseteq\Gamma\) is called \emph{nonsingular} if
\[
 \smin(A'(x))>0
 \quad\text{for Lebesgue-almost every }x\in\RT.
\]
\end{definition*}
By Lemma~\ref{lem:local}, local Lipschitz stability holds at
almost every point of \(\RT\) whenever \(T\) is nonsingular.

Let \(Q_T(x)\) be the submatrix of \(A'(x)^{\mathsf T}A'(x)\) with rows
and columns indexed by \(T\),
\[
 Q_T(x)=\bigl[A'(x)^{\mathsf T}A'(x)\bigr]_{T,T}
       \in\mathbb R^{T\times T}.
\]
This positive semidefinite matrix satisfies
\begin{equation}
 \lambda_{\min}(Q_T(x))
 =\lambda_{\min,T}\!\left(A'(x)^{\mathsf T}A'(x)\right)
 =\smin(A'(x))^2.
 \label{eq:gram}
\end{equation}
The following proposition gives a criterion for nonsingular supports.

\begin{proposition}
\label{prop:generic}
A nonempty support \(T\subseteq\Gamma\) is nonsingular if and only if
there exists \(\bar x\in\RT\) such that
\[
 \smin(A'(\bar x))>0.
\]
\end{proposition}
\begin{proof}
The existence of \(\bar x\) follows from the definition when \(T\) is
nonsingular, so it suffices to prove the converse. Let \(\bar x\in\RT\)
satisfy \(\smin(A'(\bar x))>0\). The map \(x\mapsto\det Q_T(x)\) is a polynomial in the \(|T|\) free coordinates
\((x_i)_{i\in T}\). By \eqref{eq:gram}, \(\det Q_T(\bar x)>0\), so this
polynomial is nonzero and its zero set has Lebesgue measure zero in \(\RT\).
Since \(Q_T(x)\) is positive semidefinite, \eqref{eq:gram} also gives
\[
 \{x\in\RT:\smin(A'(x))=0\}
 =\{x\in\RT:\det Q_T(x)=0\}.
\]
We conclude that \(\smin(A'(x))>0\) for almost every \(x\in\RT\), that is,
\(T\) is nonsingular.
\end{proof}

\begin{theorem}[Generic stability on a random support]
\label{thm:generic}
Fix \(d\geq1\), \(\varepsilon\in(0,1]\), and \(C>0\), and set
\[
 r=\lceil\varepsilon^{-1}\rceil,\qquad k=(2+2^d)r.
\]
Let \(T\) be uniform among the \(\binom Nt\) subsets of \(\Gamma\) of
cardinality \(t\), where \(k\leq t\leq CN^{1-\varepsilon}\).
There exist constants \(K,N_*\), depending only on \(d,\varepsilon,C\), such
that for \(N\geq N_*\), with probability at least
\(1-KN^{-\varepsilon}\), \(T\) is nonsingular.

In particular, local Lipschitz stability holds at almost every point
of \(\RT\).
\end{theorem}

The proof of Theorem~\ref{thm:generic} is presented in Appendix~\ref{app:generic}. The key idea is to construct a signal
at which the Jacobian is injective on \(\RT\). With the stated
probability, there is an \(\bar{S}\subseteq T\) such that
\[
 |\bar{S}|=k,\qquad \smin(A'(\mathbf1_{\bar{S}}))\geq1.
\]
Here \(\mathbf1_{\bar{S}}\in\mathbb R^\Gamma\) equals \(1\) on
\(\bar{S}\) and \(0\) on \(\Gamma\setminus\bar{S}\).
Proposition~\ref{prop:generic} then shows that \(T\) is nonsingular.

\section{Quantitative stability via the graph of isolated autocorrelations}
\label{sec:graph}
In this section, we introduce the graph of isolated autocorrelations for a given support set, and use it to derive explicit bounds for the Lipschitz constant in the stability estimate.  

\subsection{Graph notation and the signless Laplacian}
\label{subsec:graph-preliminaries}
\paragraph{Graph notation.} Let $G=(V,E)$ be a finite, simple, undirected graph. We denote the degree of a vertex $i$ by $d_G(i)$, the maximum degree of $G$ by $\Delta(G)$, and the complement of $G$ on $V$ by $\overline{G}$. The complete graph on $V$, in which every pair
of distinct vertices is joined by an edge, is denoted as $K_V$.

The signless Laplacian of \(G\) is denoted by \(L_+(G)=D_G+\mathcal A_G\),
where \(D_G\) is the diagonal degree matrix and \(\mathcal A_G\) denotes the
adjacency matrix.
We have, for \(u\in\mathbb R^V\), 
\begin{equation}
 u^{\mathsf T}L_+(G)u=\sum_{\{i,j\}\in E}(u_i+u_j)^2.
 \label{eq:signless}
\end{equation}
In particular $L_+(G)$ is positive semidefinite.

\begin{lemma}
\label{lem:signless}
For a graph \(G\) on \(k\geq1\) vertices,
\begin{equation}
 \lambda_{\min}(L_+(G))\geq
       \bigl(k-2-2\Delta(\overline G)\bigr)_+,
 \label{eq:signless-bound}
\end{equation}
where \(a_+=\max\{0,a\}\) for \(a\in\mathbb R\).
\end{lemma}

\begin{proof}
For \(u\in\mathbb R^V\), \eqref{eq:signless} and expansion
of squares give
\[
 u^{\mathsf T}L_+(K_V)u
 =\sum_{\substack{\{i,j\}\subseteq V\\i\neq j}}(u_i+u_j)^2
 =(k-2)\|u\|_2^2+\left(\sum_{i\in V}u_i\right)^2.
\]
For the complement \(\overline G\), we have
\[
 u^{\mathsf T}L_+(\overline G)u
 =\sum_{\substack{\{i,j\}\subseteq V\\i\neq j,\ \{i,j\}\notin E}}(u_i+u_j)^2
 \leq2\sum_{i\in V}d_{\overline G}(i)u_i^2
 \leq2\Delta(\overline G)\|u\|_2^2.
\]

Since \(G\) and \(\overline G\) partition the edges of \(K_V\),
we have \(L_+(G)=L_+(K_V)-L_+(\overline G)\). Hence
\[
 u^{\mathsf T}L_+(G)u
 =u^{\mathsf T}L_+(K_V)u
   -u^{\mathsf T}L_+(\overline G)u
 \geq\bigl(k-2-2\Delta(\overline G)\bigr)\|u\|_2^2.
\]
The left-hand side is also nonnegative by \eqref{eq:signless}.
Taking the minimum over unit vectors yields
\[
 \lambda_{\min}(L_+(G))
 =\min_{\substack{u\in\mathbb R^V\\\|u\|_2=1}}
       u^{\mathsf T}L_+(G)u
 \geq\bigl(k-2-2\Delta(\overline G)\bigr)_+.
\]

\end{proof}

\subsection{The graph of isolated autocorrelations}
\label{subsec:isolated-graph}
We introduce the graph of isolated autocorrelations, and use it to derive a
lower bound for \(\smin(A'(x))\).

\begin{definition}[Graph of isolated autocorrelations]
\label{def:isolated-graph}
For a nonempty support constraint \(T\subseteq\Gamma\), let \(E(T)\) be the
set of unordered pairs \(\{i,j\}\subseteq T\) with \(i\neq j\) such that, for
all distinct \(i_1,j_1\in T\),
\[
 i_1-j_1=\pm(i-j)
 \quad\Longrightarrow\quad
 \{i_1,j_1\}=\{i,j\}.
\]
We define the \emph{graph of isolated autocorrelations} as
\(G(T)=(T,E(T))\). We say that \(T\) is \emph{collision-free} if \(G(T)=K_T\), the complete graph on \(T\); equivalently, every pair of distinct points in \(T\) has a difference that is unique up to sign.
\end{definition}

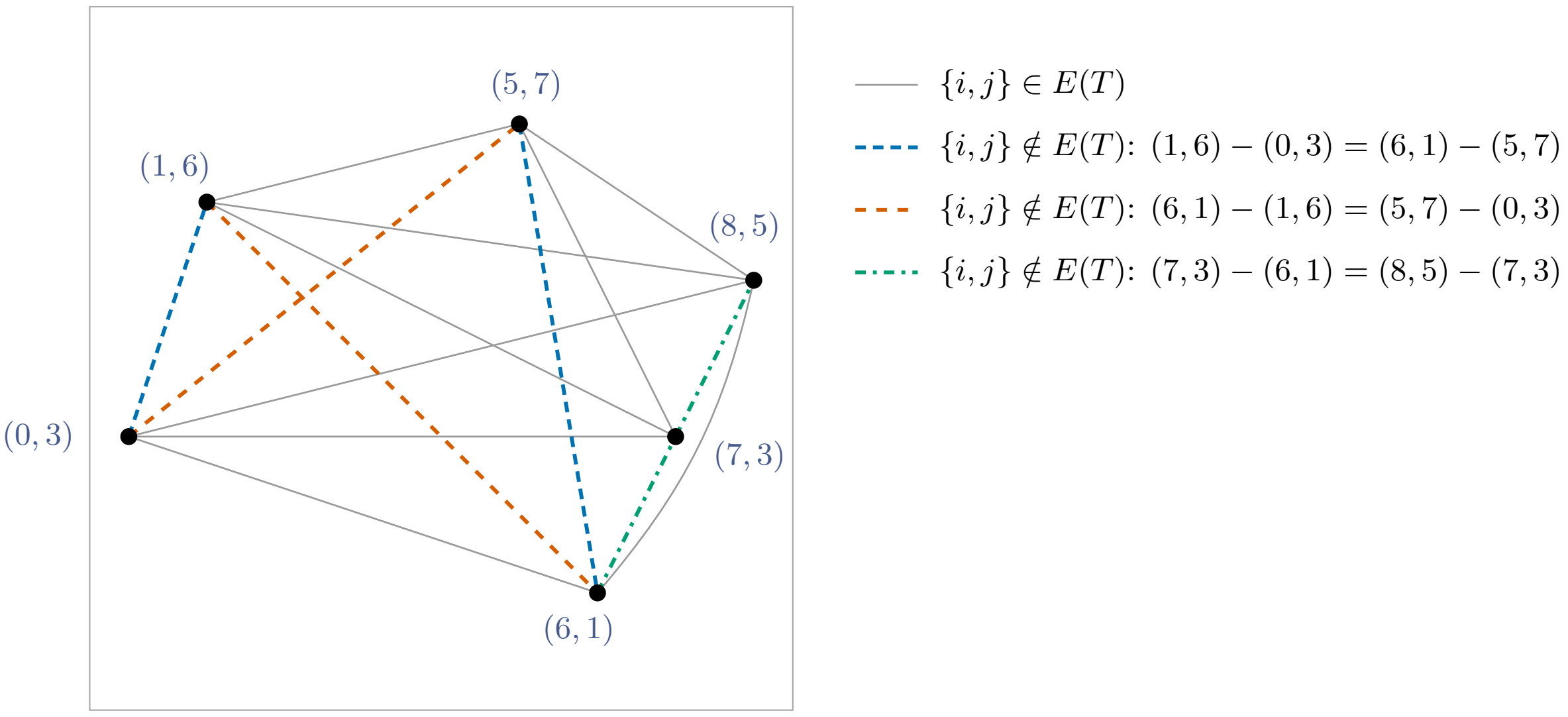
\begin{figure}[t]
 \centering
 \begin{tikzpicture}[x=0.78cm,y=0.78cm,
  pt/.style={circle,fill=black,inner sep=1.7pt},
  lab/.style={font=\footnotesize,text=gtlabel,inner sep=1.5pt},
  edge/.style={black!40,line width=0.5pt},
  miss/.style={line width=1.1pt}]
\draw[black!35,line width=0.4pt] (-0.5,-0.5) rectangle (8.5,8.5);
\coordinate (v03) at (0,3);
\coordinate (v16) at (1,6);
\coordinate (v57) at (5,7);
\coordinate (v61) at (6,1);
\coordinate (v73) at (7,3);
\coordinate (v85) at (8,5);
\draw[edge] (v03) -- (v61);
\draw[edge] (v03) -- (v73);
\draw[edge] (v03) -- (v85);
\draw[edge] (v16) -- (v57);
\draw[edge] (v16) -- (v73);
\draw[edge] (v16) -- (v85);
\draw[edge] (v57) -- (v73);
\draw[edge] (v57) -- (v85);
\draw[edge] (v61) to[bend right=14] (v85);
\draw[miss,gtblue,densely dashed] (v03) -- (v16);
\draw[miss,gtblue,densely dashed] (v57) -- (v61);
\draw[miss,gtorange,dashed] (v03) -- (v57);
\draw[miss,gtorange,dashed] (v16) -- (v61);
\draw[miss,gtgreen,dash dot] (v61) -- (v73);
\draw[miss,gtgreen,dash dot] (v73) -- (v85);
\node[pt] at (v03) {};
\node[lab,anchor=0] at ($(v03)+(180:0.62)$) {$(0,3)$};
\node[pt] at (v16) {};
\node[lab,anchor=313] at ($(v16)+(133:0.22)$) {$(1,6)$};
\node[pt] at (v57) {};
\node[lab,anchor=260] at ($(v57)+(80:0.22)$) {$(5,7)$};
\node[pt] at (v61) {};
\node[lab,anchor=63] at ($(v61)+(243:0.22)$) {$(6,1)$};
\node[pt] at (v73) {};
\node[lab,anchor=165] at ($(v73)+(345:0.42)$) {$(7,3)$};
\node[pt] at (v85) {};
\node[lab,anchor=280] at ($(v85)+(100:0.4)$) {$(8,5)$};
\draw[edge] (9.3,7.5) -- (10.100000000000001,7.5);
\node[anchor=west,font=\footnotesize] at (10.200000000000001,7.5) {$\{i,j\}\in E(T)$};
\draw[miss,gtblue,densely dashed] (9.3,6.7) -- (10.100000000000001,6.7);
\node[anchor=west,font=\footnotesize] at (10.200000000000001,6.7) {$\{i,j\}\notin E(T)$: $(1,6)-(0,3)=(6,1)-(5,7)$};
\draw[miss,gtorange,dashed] (9.3,5.9) -- (10.100000000000001,5.9);
\node[anchor=west,font=\footnotesize] at (10.200000000000001,5.9) {$\{i,j\}\notin E(T)$: $(6,1)-(1,6)=(5,7)-(0,3)$};
\draw[miss,gtgreen,dash dot] (9.3,5.1) -- (10.100000000000001,5.1);
\node[anchor=west,font=\footnotesize] at (10.200000000000001,5.1) {$\{i,j\}\notin E(T)$: $(7,3)-(6,1)=(8,5)-(7,3)$};
\end{tikzpicture}
 \caption{The graph \(G(T)\) for a support \(T\subseteq\mathbb Z_9^2\); the
 frame is \(\{0,\dots,8\}^2\). Each dashed pair has the same difference, up
 to sign, as another pair in \(T\). Differences are taken in
 \(\mathbb Z_9^2\), e.g.\ \((6,1)-(5,7)=(1,-6)=(1,3)=(1,6)-(0,3)\).}
 \label{fig:isolated-graph}
\end{figure}

Figure~\ref{fig:isolated-graph} illustrates
Definition~\ref{def:isolated-graph} for a support \(T\) of six points
in \(\Gamma=\mathbb Z_9^2\).

For \(x\in\mathbb R^\Gamma\), define
\begin{equation}
 Q_{G(T)}(x)=\sum_{\{i,j\}\in E(T)}
 (x_j\mathbf e_i+x_i\mathbf e_j)
 (x_j\mathbf e_i+x_i\mathbf e_j)^{\mathsf T}
 \in\mathbb R^{\Gamma\times\Gamma},
 \label{eq:graph-matrix}
\end{equation}
where \(\mathbf e_i\in\mathbb R^\Gamma\) is the coordinate vector at \(i\).
This matrix is positive semidefinite, with
\begin{equation}
 h^{\mathsf T}Q_{G(T)}(x)h=\sum_{\{i,j\}\in E(T)}(x_jh_i+x_ih_j)^2.
 \label{eq:graph-quadratic-form}
\end{equation}

\begin{lemma}
\label{lem:graph-comparison}
For \(x\in\RT\), the matrix \(Q_{G(T)}(x)\) in \eqref{eq:graph-matrix}
satisfies
\begin{equation}
 \smin(A'(x))^2\geq2\lambda_{\min,T}(Q_{G(T)}(x)).
 \label{eq:comparison}
\end{equation}
\end{lemma}
\begin{proof}
We first show that
\begin{equation}
 (A'(x)h)_{v_1}
 =\begin{cases}
 x_{j_1}h_{i_1}+x_{i_1}h_{j_1},&v_1\neq-v_1,\\
 2(x_{j_1}h_{i_1}+x_{i_1}h_{j_1}),&v_1=-v_1
 \end{cases}
 \label{eq:isolated-rows}
\end{equation}
for all \(h\in\RT\), \(\{i_1,j_1\}\in E(T)\), and
\(v_1\in\{i_1-j_1,\,j_1-i_1\}\). By Definition~\ref{def:isolated-graph},
\[
 \{(i,j):\ i,j\in T,\ i-j=v_1\}\subseteq\{(i_1,j_1),(j_1,i_1)\}.
\]
Therefore \eqref{eq:Ah-pairs} gives\footnote{Here \(\mathbf1_{\{\cdot\}}\) equals \(1\)
if the condition in braces holds and \(0\) otherwise.}
\begin{align*}
 (A'(x)h)_{v_1}
 &=\sum_{\substack{i,j\in T\\i-j=v_1}}(x_jh_i+x_ih_j)\\
 &=\mathbf1_{\{i_1-j_1=v_1\}}\,(x_{j_1}h_{i_1}+x_{i_1}h_{j_1})
  +\mathbf1_{\{j_1-i_1=v_1\}}\,(x_{i_1}h_{j_1}+x_{j_1}h_{i_1})\\
 &=\bigl(\mathbf1_{\{i_1-j_1=v_1\}}+\mathbf1_{\{j_1-i_1=v_1\}}\bigr)
   (x_{j_1}h_{i_1}+x_{i_1}h_{j_1}).
\end{align*}
This gives \eqref{eq:isolated-rows}.

Next, by Definition~\ref{def:isolated-graph}, for
\(\{i,j\},\{i',j'\}\in E(T)\) with \(\{i,j\}\neq\{i',j'\}\), we have
\[
 \{i-j,\,j-i\}\cap\{i'-j',\,j'-i'\}=\varnothing.
\]
Combined with \(\bigcup_{\{i,j\}\in E(T)}\{i-j,\,j-i\}\subseteq\Gamma\),
this gives, for any \(h\in\RT\),
\begin{align*}
 h^{\mathsf T}A'(x)^{\mathsf T}A'(x)h
 &=\sum_{v\in\Gamma}(A'(x)h)_v^2
  \geq\sum_{\{i,j\}\in E(T)}\sum_{v\in\{i-j,\,j-i\}}(A'(x)h)_v^2\\
 &=\sum_{\substack{\{i,j\}\in E(T)\\i-j\neq j-i}}
    \bigl[(A'(x)h)_{i-j}^2+(A'(x)h)_{j-i}^2\bigr]
   +\sum_{\substack{\{i,j\}\in E(T)\\i-j=j-i}}(A'(x)h)_{i-j}^2\\
 &=2\sum_{\substack{\{i,j\}\in E(T)\\i-j\neq j-i}}(x_jh_i+x_ih_j)^2
   +4\sum_{\substack{\{i,j\}\in E(T)\\i-j=j-i}}(x_jh_i+x_ih_j)^2\\
 &\geq2h^{\mathsf T}Q_{G(T)}(x)h.
\end{align*}
Taking the infimum over unit \(h\in\RT\) gives \eqref{eq:comparison}.
\end{proof}

\subsection{Quantitative local stability}

In this subsection, we combine Lemmas~\ref{lem:signless}
and~\ref{lem:graph-comparison} to give an explicit lower bound for
\(\smin(A'(x^*))\) at every nonzero \(x^*\in\RT\).

\begin{theorem}[Quantitative local stability]
\label{thm:quantitative}
Let \(x^*\neq0\), \(S=\supp(x^*)\subseteq T\), and \(s=|S|\). Set
\[
 \Delta=\Delta(\overline{G(T)}),\qquad
 m=\min_{i\in S}|x_i^*|,\quad
 M=\max_{i\in S}|x_i^*|.
\]
Then
\begin{equation}
 \smin(A'(x^*))
 \geq\sqrt2\,\frac{m^2}{M}\sqrt{(s-2-2\Delta)_+}.
 \label{eq:quantitative}
\end{equation}
If \(s>2+2\Delta\), there exists a convex neighborhood \(\mathcal U\) of \(x^*\) in
\(\RT\) such that
\begin{equation}
 \|x^{(1)}-x^{(2)}\|_2\leq
 \frac{\sqrt2\,M}{m^2}\sqrt{\frac{N}{s-2-2\Delta}}\,
 \bigl\|\,|\mathcal F(x^{(1)})|^2-|\mathcal F(x^{(2)})|^2\,\bigr\|_2,
 \qquad x^{(1)},x^{(2)}\in\mathcal U.
 \label{eq:quantitative-local}
\end{equation}
\end{theorem}
\begin{proof}
We fix \(h\in\RT\) and write \(Z=T\setminus S\). Since \(x^*\) vanishes on
\(Z\), \eqref{eq:graph-quadratic-form} gives
\begin{equation}
 \begin{aligned}
 &h^{\mathsf T}Q_{G(T)}(x^*)h
 =\sum_{\{i,j\}\in E(T)}(x_j^*h_i+x_i^*h_j)^2\\
 &=\sum_{\substack{\{i,j\}\in E(T)\\i,j\in S}}(x_j^*h_i+x_i^*h_j)^2
 +\sum_{\substack{i\in S,\ j\in Z\\\{i,j\}\in E(T)}}
       (\underbrace{x_j^*}_{=0}h_i+x_i^*h_j)^2
 +\sum_{\substack{\{i,j\}\in E(T)\\i,j\in Z}}
       (\underbrace{x_j^*}_{=0}h_i+\underbrace{x_i^*}_{=0}h_j)^2\\
 &=\underbrace{\sum_{\substack{\{i,j\}\in E(T)\\i,j\in S}}
       (x_j^*h_i+x_i^*h_j)^2}_{\text{(i)}}
 +\underbrace{\sum_{j\in Z}h_j^2
       \sum_{\substack{i\in S\\\{i,j\}\in E(T)}}(x_i^*)^2}_{\text{(ii)}}.
 \end{aligned}
 \label{eq:split}
\end{equation}

We first bound (i). Let
\[
 G_S=\bigl(S,\ \bigl\{\{i,j\}\in E(T):\ i,j\in S\bigr\}\bigr)
\]
be the subgraph of \(G(T)\) induced by \(S\), and let
\(u=\bigl(h_i/x_i^*\bigr)_{i\in S}\in\mathbb R^S\). We write (i) as a
signless Laplacian quadratic form:
\[
 \begin{aligned}
 \text{(i)}
 &=\sum_{\{i,j\}\in E(G_S)}(x_j^*h_i+x_i^*h_j)^2
  =\sum_{\{i,j\}\in E(G_S)}(x_i^*x_j^*)^2
   \Bigl(\frac{h_i}{x_i^*}+\frac{h_j}{x_j^*}\Bigr)^2\\
 &=\sum_{\{i,j\}\in E(G_S)}(x_i^*x_j^*)^2(u_i+u_j)^2
  \geq m^4\sum_{\{i,j\}\in E(G_S)}(u_i+u_j)^2
  \xlongequal{\eqref{eq:signless}}m^4\,u^{\mathsf T}L_+(G_S)u,
 \end{aligned}
\]
where the inequality above uses \(|x_i^*|\geq m\) on
\(S\).

Next, for \(i\in S\), every vertex of \(S\) not adjacent to \(i\) in
\(G_S\) is not adjacent to \(i\) in \(G(T)\), so
\(\Delta(\overline{G_S})\leq\Delta(\overline{G(T)})=\Delta\), and
Lemma~\ref{lem:signless} with \(k=s\) gives
\[
 m^4\,u^{\mathsf T}L_+(G_S)u
 \geq m^4(s-2-2\Delta)_+\sum_{i\in S}u_i^2
 =m^4(s-2-2\Delta)_+\sum_{i\in S}\frac{h_i^2}{(x_i^*)^2}
 \geq\frac{m^4}{M^2}(s-2-2\Delta)_+\sum_{i\in S}h_i^2,
\]
where the last inequality uses \(|x_i^*|\leq M\) on \(S\).

We then bound (ii). Since \(d_{\overline{G(T)}}(j)\leq\Delta\), each
\(j\in Z\) is adjacent in \(G(T)\) to at least \(s-\Delta\) points of
\(S\). Thus, in \eqref{eq:split}, the inner sum of (ii) has at least
\(s-\Delta\) terms, each at least \(m^2\) since \(|x_i^*|\geq m\) on
\(S\). Hence
\[
 \text{(ii)}
 \geq m^2(s-\Delta)_+\sum_{j\in Z}h_j^2
 \geq\frac{m^4}{M^2}(s-2-2\Delta)_+\sum_{j\in Z}h_j^2.
\]

Finally, adding the bounds on (i) and (ii), we obtain
\[
 h^{\mathsf T}Q_{G(T)}(x^*)h
 =\text{(i)}+\text{(ii)}
 \geq\frac{m^4}{M^2}(s-2-2\Delta)_+
 \Bigl(\sum_{i\in S}h_i^2+\sum_{j\in Z}h_j^2\Bigr)
 =\frac{m^4}{M^2}(s-2-2\Delta)_+\|h\|_2^2.
\]
Since \(h\in\RT\) is arbitrary, Lemma~\ref{lem:graph-comparison} gives
\[
 \smin(A'(x^*))
 \geq\sqrt{2\lambda_{\min,T}\bigl(Q_{G(T)}(x^*)\bigr)}
 \geq\sqrt2\,\frac{m^2}{M}\sqrt{(s-2-2\Delta)_+}.
\]

If \(s>2+2\Delta\), then \eqref{eq:quantitative} gives
\(\sigma=\smin(A'(x^*))\geq\sqrt2\,\frac{m^2}{M}\sqrt{s-2-2\Delta}>0\), and
applying Lemma~\ref{lem:local} with \(L=\frac2\sigma\) and
\eqref{eq:normalization}, we obtain \eqref{eq:quantitative-local}.
\end{proof}

\begin{remark}
\label{rem:collisionfree}
If \(T\) is collision-free, then \(G(T)=K_T\), so \(\Delta=0\), and
Theorem~\ref{thm:quantitative} gives, for \(s=|\supp(x^*)|\geq3\),
\[
 \smin(A'(x^*))\geq\sqrt{2(s-2)}\,\frac{m^2}{M}.
\]
For fixed \(m\) and \(M\), this bound has the same \(\sqrt s\) scaling as
the collision-free estimate in \cite[Lemma~5]{ghosh2023sparse}.
\end{remark}

\subsection{Nonnegativity implies global stability}
\label{subsec:positive-global}

We now restrict to nonnegative signals \(x\in\RT\).
To obtain a global estimate, fix an integer \(1\leq s\leq|T|\)
and \(0<m\leq M\), and define
\[
 \mathcal K=\left\{x\in\RT\cap\mathbb R^\Gamma_+:
 |\supp(x)|\geq s,\quad
 m\leq x_i\leq M\ \text{for }i\in\supp(x)\right\}.
\]

\begin{theorem}[Global stability under nonnegativity]
\label{thm:positive-global}
Let \(\Delta=\Delta(\overline{G(T)})\). If
\(s>2+2\Delta\), then for all \(x^{(1)},x^{(2)}\in\mathcal K\),
\begin{equation}
 \|x^{(1)}-x^{(2)}\|_2
 \leq\frac{2\sqrt2\,M}{m^2}\sqrt{\frac{N}{s-2-2\Delta}}\,
 \bigl\|\,|\mathcal F(x^{(1)})|^2-|\mathcal F(x^{(2)})|^2\,\bigr\|_2.
 \label{eq:positive-global}
\end{equation}
\end{theorem}
\begin{proof}
We fix \(x^{(1)},x^{(2)}\in\mathcal K\) and write
\(S_j=\supp(x^{(j)})\) for \(j=1,2\). By nonnegativity, the midpoint has support
\(S_1\cup S_2\subseteq T\), with
\[
 |S_1\cup S_2|\geq s,\qquad
 \frac m2\leq\frac{x^{(1)}_i+x^{(2)}_i}{2}\leq M
 \quad\text{for every } i\in S_1\cup S_2.
\]
Applying Theorem~\ref{thm:quantitative} at \((x^{(1)}+x^{(2)})/2\), with
the same support constraint \(T\), and using
\(|S_1\cup S_2|-2-2\Delta\geq s-2-2\Delta>0\), we obtain
\[
 \smin\Bigl(A'\Bigl(\frac{x^{(1)}+x^{(2)}}{2}\Bigr)\Bigr)
 \geq\sqrt2\,\frac{(m/2)^2}{M}
       \sqrt{s-2-2\Delta}
 =\frac{m^2\sqrt{s-2-2\Delta}}{2\sqrt2\,M}.
\]
Thus Lemma~\ref{lem:global} applies with
\(\sigma=\frac{m^2\sqrt{s-2-2\Delta}}{2\sqrt2\,M}\), and
\eqref{eq:global-single-sigma} and \eqref{eq:normalization} give
\eqref{eq:positive-global}.
\end{proof}

\section{Quantitative stability for random supports}
\label{sec:random}

\begin{samepage}
In this section, \(T\subseteq\Gamma\) is a Bernoulli random support:
each point of \(\Gamma\) is included in \(T\) independently with
probability \(\frac tN\), so \(\mathbb E|T|=t\). We first show that, for
\(t=O(\sqrt N)\), with high probability \(|T|\) is close
to \(t\) and \(\Delta(\overline{G(T)})\leq\theta t\). On this event, we
then derive local stability bounds at every \(x\in\RT\) with
\(|\supp(x)|=\Omega(|T|)\), and global stability bounds for such \(x\)
that are nonnegative.

\subsection{Bounds for \texorpdfstring{\(|T|\) and
\(\Delta(\overline{G(T)})\)}{|T| and the maximum degree of the complement of G(T)}}
\end{samepage}

\begin{lemma}[Bounds for \(|T|\) and \(\Delta(\overline{G(T)})\)]
\label{lem:random-degree}
Let \(T\subseteq\Gamma\) be a Bernoulli random support with
\(\mathbb E|T|=t\).
Fix \(0<\theta<\frac12\), and assume
\[
 (\log N)^4=o(t),\qquad t\leq c\sqrt N,
\]
where \(c\) is a constant satisfying
\(0<c<\sqrt{\frac{\theta}{6}}\).
There exist positive constants \(C_1,c_1,C_2,c_2\), depending only on
\(d,c,\theta\), such that, for all sufficiently large \(N\), with probability
at least \(1-\varepsilon_\theta(t,N)\),
\begin{equation}
 (1-\theta)t\leq|T|\leq(1+\theta)t,
 \qquad \Delta(\overline{G(T)})\leq\theta t,
 \label{eq:random-degree}
\end{equation}
where
\begin{equation}
 \varepsilon_\theta(t,N)=2e^{-\theta^2t/3}
       +C_1t e^{-c_1\sqrt t\log(N/t)}+C_2t e^{-c_2t}.
 \label{eq:failure}
\end{equation}
\end{lemma}

The proof is given in Appendix~\ref{app:degree}.
Since \((\log N)^4=o(t)\) forces \(t\to\infty\), and \(t\leq c\sqrt N\)
gives \(\log(N/t)\geq\frac12\log N\), each term of \eqref{eq:failure}
tends to zero.

\subsection{Local stability for random supports}

For fixed \(\alpha\in(0,1]\) and \(0<m\leq M\), we study
local stability at points in the class
\[
 \mathcal X_{T,\alpha,m,M}=\left\{x\in\RT:
 |\supp(x)|\geq\alpha|T|,\quad
 m\leq|x_i|\leq M\ \text{for }i\in\supp(x)
 \right\}.
\]

\begin{theorem}
\label{thm:random-stability}
Let \(\alpha,m,M\) be as above, set \(\theta=\frac{\alpha}{2(2+\alpha)}\),
and fix \(0<c<\sqrt{\frac{\theta}{6}}\). Let \(T\subseteq\Gamma\) be a Bernoulli
random support with \(\mathbb E|T|=t\), where
\[
 (\log N)^4=o(t),\qquad t\leq c\sqrt N.
\]
Then, for all sufficiently large \(N\), with probability at least
\(1-\varepsilon_\theta(t,N)\), where \(\varepsilon_\theta\) is given by
\eqref{eq:failure}, every
\(x^*\in\mathcal X_{T,\alpha,m,M}\) satisfies
\begin{equation}
 \smin(A'(x^*))\geq\frac{m^2}{M}\sqrt{\frac{\alpha t}{2}},
 \label{eq:random-sigma}
\end{equation}
and has a convex neighborhood \(\mathcal U\) in \(\RT\) such that
\begin{equation}
 \|x^{(1)}-x^{(2)}\|_2\leq\frac{2\sqrt2\,M}{m^2}\sqrt{\frac{N}{\alpha t}}\,
 \bigl\|\,|\mathcal F(x^{(1)})|^2-|\mathcal F(x^{(2)})|^2\,\bigr\|_2,
 \qquad x^{(1)},x^{(2)}\in\mathcal U.
 \label{eq:random-local}
\end{equation}
\end{theorem}
\begin{proof}
Since \(\alpha\leq1\), we have \(\theta\leq\frac16<\frac12\), so
Lemma~\ref{lem:random-degree} applies: for all sufficiently large \(N\),
with probability at least \(1-\varepsilon_\theta(t,N)\),
\[
 |T|\geq(1-\theta)t,\qquad
 \Delta=\Delta(\overline{G(T)})\leq\theta t.
\]
It therefore suffices to show that, whenever \(T\) satisfies these two
bounds, every \(x^*\in\mathcal X_{T,\alpha,m,M}\) satisfies
\eqref{eq:random-sigma} and \eqref{eq:random-local}.

For \(x^*\in\mathcal X_{T,\alpha,m,M}\), we write \(S=\supp(x^*)\). The
bound on \(|T|\) gives
\[
 s=|S|\geq\alpha|T|
          \geq\alpha(1-\theta)t.
\]
Combining this with the bound on \(\Delta\) and
\(\theta=\frac{\alpha}{2(2+\alpha)}\), we obtain, for all sufficiently
large \(N\),
\[
 s-2-2\Delta
 \geq\bigl(\alpha(1-\theta)-2\theta\bigr)t-2
 =\frac{\alpha t}{2}-2\geq\frac{\alpha t}{4}.
\]

Finally, substituting \(m\leq|x_i^*|\leq M\) for \(i\in S\) and
\(s-2-2\Delta\geq\frac{\alpha t}{4}>0\) into
\eqref{eq:quantitative} and \eqref{eq:quantitative-local} of
Theorem~\ref{thm:quantitative} gives \eqref{eq:random-sigma} and
\eqref{eq:random-local}.
\end{proof}

\subsection{Global stability for nonnegative signals}
\label{subsec:random-positive-global}

Combining Theorem~\ref{thm:random-stability} with Lemma~\ref{lem:global},
we obtain a global stability estimate for nonnegative signals in
\(\mathcal X_{T,\alpha,m,M}\).

\begin{theorem}[Global stability on a random support]
\label{thm:random-positive-global}
Under the hypotheses and notation of Theorem~\ref{thm:random-stability},
let \(\mathcal K=\mathcal X_{T,\alpha,m,M}\cap\mathbb R^\Gamma_+\).
For all sufficiently large \(N\), with probability at least
\(1-\varepsilon_\theta(t,N)\), the estimate
\begin{equation}
 \|x^{(1)}-x^{(2)}\|_2
 \leq\frac{4\sqrt2\,M}{m^2}\sqrt{\frac{N}{\alpha t}}\,
 \bigl\|\,|\mathcal F(x^{(1)})|^2-|\mathcal F(x^{(2)})|^2\,\bigr\|_2
 \label{eq:random-positive-global}
\end{equation}
holds for all
\(x^{(1)},x^{(2)}\in\mathcal K\).
\end{theorem}
\begin{proof}
We first show that, for all \(x^{(1)},x^{(2)}\in\mathcal K\), the
midpoint lies in \(\mathcal X_{T,\alpha,m/2,M}\). We write
\(S_j=\supp(x^{(j)})\) for \(j=1,2\). Nonnegativity gives
\[
 \supp\Bigl(\frac{x^{(1)}+x^{(2)}}{2}\Bigr)=S_1\cup S_2,
 \qquad\text{so}\qquad
 |S_1\cup S_2|\geq|S_1|\geq\alpha|T|.
\]
Moreover, \(m/2\leq(x^{(1)}_i+x^{(2)}_i)/2\leq M\) for every
\(i\in S_1\cup S_2\).
Thus \((x^{(1)}+x^{(2)})/2\in\mathcal X_{T,\alpha,m/2,M}\).

We next apply Theorem~\ref{thm:random-stability} with \(m/2\) in place of
\(m\). Since \(\theta\) and \(\varepsilon_\theta\) depend only on
\(\alpha\), we obtain, with probability at least
\(1-\varepsilon_\theta(t,N)\), for all
\(x^{(1)},x^{(2)}\in\mathcal K\),
\[
 \smin\Bigl(A'\Bigl(\frac{x^{(1)}+x^{(2)}}{2}\Bigr)\Bigr)
 \geq\frac{(m/2)^2}{M}\sqrt{\frac{\alpha t}{2}}
 =\frac{m^2\sqrt{\alpha t}}{4\sqrt2\,M}.
\]
Thus Lemma~\ref{lem:global} applies with
\(\sigma=\frac{m^2\sqrt{\alpha t}}{4\sqrt2\,M}\), and
\eqref{eq:global-single-sigma} and \eqref{eq:normalization} give
\eqref{eq:random-positive-global}.
\end{proof}

\section{Numerical experiments}
\label{sec:numerics}

We present two numerical experiments.
In Subsection~\ref{subsec:numerics-graph}, we investigate whether our
quantitative lower bound on \(\smin(A'(x))\) remains positive for
sparsities beyond the order of \(\sqrt N\).
In Subsection~\ref{subsec:numerics-threshold}, we estimate the probability
that a uniformly random \(t\)-element subset of \(\Gamma\) is nonsingular,
for \(t\) up to the rank bound \eqref{eq:row-bound}.

\subsection{Sparsity limit of the quantitative lower bound}
\label{subsec:numerics-graph}

Recall the lower bound \eqref{eq:quantitative}: for every
\(x^*\in\RT\setminus\{0\}\),
\[
 \smin(A'(x^*))\geq\sqrt2\,\frac{m^2}{M}\sqrt{(s-2-2\Delta)_+},
\]
where \(s=|\supp(x^*)|\), \(m=\min_{i\in\supp(x^*)}|x^*_i|\),
\(M=\max_{i\in\supp(x^*)}|x^*_i|\), and \(\Delta=\Delta(\overline{G(T)})\).
Since \(s\leq|T|\), a necessary condition for this lower bound to be
positive at some \(x^*\in\RT\setminus\{0\}\) is
\begin{equation}
 |T|>2+2\Delta(\overline{G(T)}).
 \label{eq:numerics-positive}
\end{equation}

In this experiment, \(T\) is a Bernoulli random support: each
point of \(\Gamma\) lies in \(T\) independently with probability
\(t/N\), so that \(t=\mathbb E|T|\). We use \(N=n^2\) with
\(n\in\{100,300,900\}\), both for \(\Gamma=\mathbb Z_N\) and for
\(\Gamma=\mathbb Z_n^2\). For each \(N\) and each
\[
 \frac t{\sqrt N}\in\{0.1,0.2,\dots,2.0\},
\]
we sample \(T\) independently \(200\) times and take the proportion of
samples satisfying \eqref{eq:numerics-positive} as its empirical
probability.

\begin{figure}[t]
 \centering
 \includegraphics[width=\linewidth]{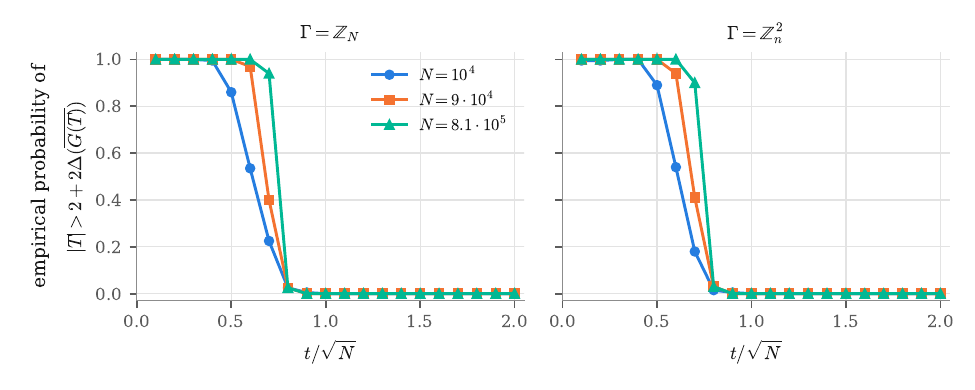}
 \caption{Empirical probability of \eqref{eq:numerics-positive}, estimated
 from \(200\) independent samples of a Bernoulli support \(T\) with
 \(\mathbb E|T|=t\), against \(t/\sqrt N\). The curves correspond to
 \(N=n^2\) with \(n=100,300,900\).}
 \label{fig:numerics-graph}
\end{figure}

Figure~\ref{fig:numerics-graph} shows the results. For
\(t/\sqrt N\leq0.5\), the empirical probability of
\eqref{eq:numerics-positive} is at least \(0.86\) for every \(N\) and for
both \(\Gamma=\mathbb Z_N\) and \(\Gamma=\mathbb Z_n^2\). For \(t/\sqrt N\geq1\), it is \(0\): none of the
\(200\) samples satisfies \eqref{eq:numerics-positive}. In between, the
probability drops, more steeply for larger \(N\), and the drop moves
slightly toward larger \(t/\sqrt N\) as \(N\) grows.

Thus, in our experiments, the bound \eqref{eq:quantitative} vanished for
every \(x^*\in\RT\setminus\{0\}\) once \(t\geq\sqrt N\). This suggests
that, up to a constant factor, \(\sqrt N\) is the limit of this
quantitative lower bound on \(\smin(A'(x^*))\). It does not mean that local stability
itself fails for larger \(t\); see
Subsection~\ref{subsec:numerics-threshold}.

\subsection{Nonsingular supports up to the rank bound}
\label{subsec:numerics-threshold}

Recall from \eqref{eq:row-bound} that \(A'(x)\) has at most
\(R=(N+|\Gamma[2]|)/2\) distinct rows; hence, by
Remark~\ref{rem:support-rank}, every support with \(|T|>R\) is singular.
In Section~\ref{sec:generic}, we proved that, for
\(t\leq CN^{1-\varepsilon}\), a uniformly random \(t\)-element subset
\(T\) of \(\Gamma\) is nonsingular with high probability. We now
investigate whether this remains true for every \(t\leq R\).

In this experiment, \(T\) is a uniformly random \(t\)-element subset of
\(\Gamma\). We use \(N=n^2\) with \(n\in\{20,30,40\}\), both for
\(\Gamma=\mathbb Z_N\) and for \(\Gamma=\mathbb Z_n^2\), and
\[
 t\in\bigl\{\lceil R/10\rceil,\lceil2R/10\rceil,\dots,\lceil11R/10\rceil\bigr\}.
\]
For each \((N,t)\), we sample \(T\) independently \(200\) times.

For each sampled \(T\), let \(x^*=g/\|g\|_2\), where
\[
 g_i\stackrel{\mathrm{iid}}{\sim}\mathcal N(0,1)\quad(i\in T),\qquad
 g_i=0\quad(i\in\Gamma\setminus T).
\]
Then \(x^*\) is a unit vector in \(\RT\).
We compute \(\smin(A'(x^*))\) by a singular value decomposition. By
Proposition~\ref{prop:generic}, \(T\) is nonsingular whenever
\(\smin(A'(x^*))>0\). We count \(T\) as nonsingular if the computed value
exceeds \(10^{-8}\).
The largest size \(t=\lceil11R/10\rceil>R\) serves as a comparison: there
\(\smin(A'(x^*))=0\) by Remark~\ref{rem:support-rank}, so its computed
values show the size of rounding errors.

\begin{figure}[t]
 \centering
 \includegraphics[width=\linewidth]{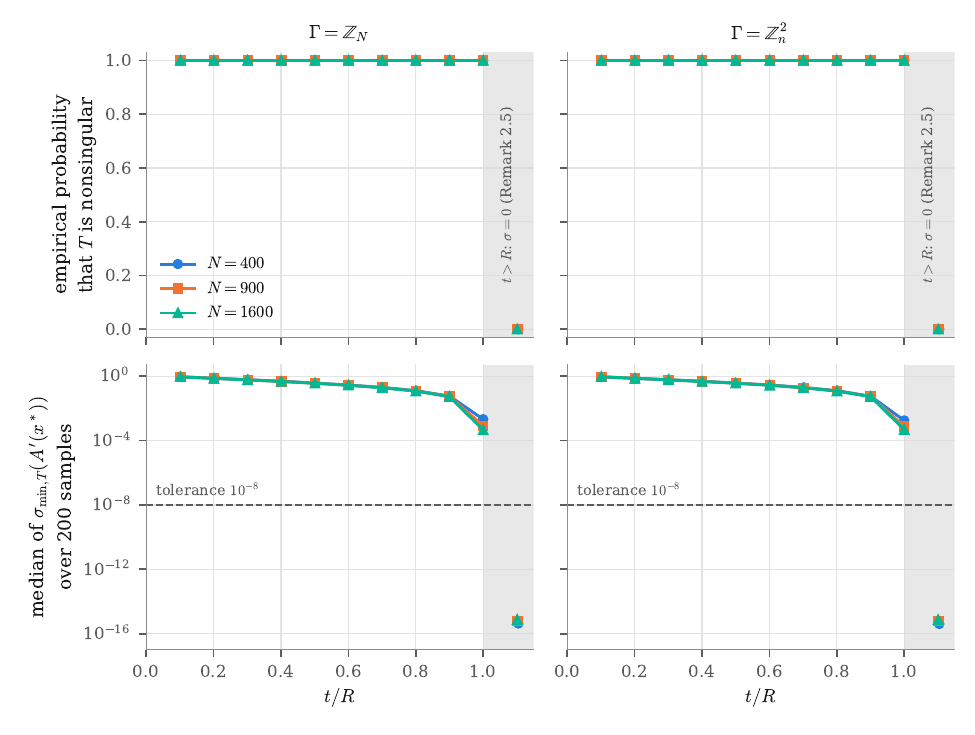}
 \caption{Uniformly random \(t\)-element subsets \(T\) of \(\Gamma\), with
 \(200\) samples for each \(t\), against \(t/R\). The curves correspond to
 \(N=n^2\) with \(n=20,30,40\). Upper panels: empirical probability that
 \(T\) is nonsingular. Lower panels: median of \(\smin(A'(x^*))\) over
 the samples; the dashed line is the tolerance \(10^{-8}\). In the shaded
 region \(t>R\), no subset is nonsingular by
 Remark~\ref{rem:support-rank}; the points there are the comparison case
 \(t=\lceil11R/10\rceil\).}
 \label{fig:numerics-threshold}
\end{figure}

For every \(N\), for both \(\Gamma=\mathbb Z_N\) and
\(\Gamma=\mathbb Z_n^2\), and for every \(t\leq R\), all \(200\) samples
of \(T\) were nonsingular (Figure~\ref{fig:numerics-threshold}). For
\(t=\lceil11R/10\rceil\), all computed values were below \(10^{-15}\), far
below the tolerance \(10^{-8}\).

Based on these results, we conjecture that the necessary condition
\(|T|\leq R\) of Remark~\ref{rem:support-rank} is also sufficient: for any
\(t=t(N)\leq R\), a uniformly random \(t\)-element subset of \(\Gamma\) is
nonsingular with probability tending to one as \(N\to\infty\). This would
extend Theorem~\ref{thm:generic} from \(t\leq CN^{1-\varepsilon}\) to all
\(t\leq R\).

Since \(|\Gamma[2]|\leq2^d\), \(R/N\to\frac12\) as \(N\to\infty\), so
the conjectured threshold \(R\) is asymptotically \(\frac N2\).

\section*{Declaration of generative AI and AI-assisted technologies in the manuscript preparation process}
We used Claude Opus 5.5 to improve the language and grammar of
the manuscript. Claude Opus 5.5 also helped write the code for the
numerical experiments in Section~\ref{sec:numerics}. GPT-6 Astra assisted
with the proof of
Lemma~\ref{lem:random-degree}, in particular in connecting our probability
estimates to Warnke's upper-tail bound~\cite{Warnke2017UpperTails}. 

\clearpage
\appendix
\section{Proof of Theorem~\ref{thm:generic}}
\label{app:generic}

By Proposition~\ref{prop:generic}, \(T\) is nonsingular as soon as
\(\smin(A'(\bar x))>0\) at a single point \(\bar x\in\RT\). We take
\(\bar x=\mathbf1_{\bar{S}}\) for a \(k\)-element subset
\(\bar{S}\subseteq T\), and write \(\bar{Z}=T\setminus\bar{S}\).

We generate the pair \((\bar{S},\bar{Z})\) as follows: \(\bar{S}\) is a
uniformly random \(k\)-element subset of \(\Gamma\), and, given
\(\bar{S}\), \(\bar{Z}\) is a uniformly random \((t-k)\)-element subset of
\(\Gamma\setminus\bar{S}\). We write \(\mathbb P\) for probability under
this sampling. For disjoint \(S_0,Z_0\subseteq\Gamma\) with \(|S_0|=k\) and
\(|Z_0|=t-k\),
\begin{equation}
 \mathbb P(\bar{S}=S_0,\bar{Z}=Z_0)
 =\frac{1}{\binom Nk\binom{N-k}{t-k}}
 =\frac{1}{\binom Nt\binom tk},
 \label{eq:same-law}
\end{equation}
which is the probability of first choosing a uniformly random
\(t\)-element set \(T=S_0\cup Z_0\) and then \(\bar{S}=S_0\) inside it.
Hence \(T=\bar{S}\cup\bar{Z}\) is a uniformly random \(t\)-element subset
of \(\Gamma\), as in Theorem~\ref{thm:generic}, and it suffices to show
that \(\mathbb P(\bar{S}\cup\bar{Z}\text{ is nonsingular})\geq
1-KN^{-\varepsilon}\).

The proof has three steps. In Appendix~\ref{app:construction}, we give
three deterministic conditions on \((\bar{S},\bar{Z})\) that imply
\(\smin(A'(\mathbf1_{\bar{S}}))\geq1\). In
Appendix~\ref{app:probabilities}, we bound the probability that each
condition fails. In Appendix~\ref{app:generic-proof}, we combine these
bounds to complete the proof.

\subsection{Sufficient conditions for nonsingularity of
\texorpdfstring{\(T=\bar{S}\cup\bar{Z}\)}{T = S-bar union Z-bar}}
\label{app:construction}

We now give three conditions on \((\bar{S},\bar{Z})\) that guarantee
\(\smin(A'(\mathbf1_{\bar{S}}))>0\).
For \(A,B\subseteq\Gamma\) and \(z\in\Gamma\), we write
\[
 A-B=\{a-b:\ a\in A,\ b\in B\},\quad
 z-B=\{z\}-B,\quad B-z=B-\{z\},\quad
 -A=\{-a:\ a\in A\}.
\]

\begin{lemma}[Sufficient conditions for nonsingularity of \(T=\bar{S}\cup\bar{Z}\)]
\label{lem:construct}
Let \(\bar{S},\bar{Z}\subseteq\Gamma\) be disjoint with
\(|\bar{S}|=k\geq3\), and let \(T=\bar{S}\cup\bar{Z}\). Suppose that \((\bar{S},\bar{Z})\)
satisfies the following conditions:
\begin{enumerate}
\item[(i)] The set \(\bar{S}\) is collision-free.
\par\noindent
For distinct \(a,b\in \bar{S}\)
and distinct \(a_1,b_1\in \bar{S}\), \(a-b=\pm(a_1-b_1)\) implies
\(\{a,b\}=\{a_1,b_1\}\).

\item[(ii)] \((\bar{S}-\bar{S})\cap(\bar{Z}-\bar{S})=\varnothing\).

\item[(iii)] For every \(z\in \bar{Z}\),
\begin{equation}
 z-\bar{S}\nsubseteq
 \bigcup_{w\in \bar{Z}\setminus\{z\}}
 \bigl((w-\bar{S})\cup(\bar{S}-w)\bigr).
 \label{eq:cond-iii}
\end{equation}
\end{enumerate}
Then
\[
 \smin(A'(\mathbf1_{\bar{S}}))\geq1.
\]
In particular, \(T=\bar{S}\cup\bar{Z}\) is nonsingular.
\end{lemma}
\begin{proof}
We prove that
\[
 \|A'(\mathbf1_{\bar{S}})h\|_2^2\geq\|h\|_2^2,
 \qquad h\in\RT,
\]
which gives \(\smin(A'(\mathbf1_{\bar{S}}))\geq1\).

\paragraph{Step 1.}
We apply \eqref{eq:Ah-pairs} with \(x=\mathbf1_{\bar{S}}\in\RT\) and
split \(T=\bar{S}\cup\bar{Z}\):
\begin{align}
 (A'(\mathbf1_{\bar{S}})h)_v
 &=\sum_{\substack{i,j\in T\\i-j=v}}
   \bigl((\mathbf1_{\bar{S}})_jh_i+(\mathbf1_{\bar{S}})_ih_j\bigr)
  =\sum_{\substack{i\in T,\ j\in\bar{S}\\i-j=v}}h_i
   +\sum_{\substack{i\in\bar{S},\ j\in T\\i-j=v}}h_j\notag\\
 &=\sum_{\substack{i,j\in\bar{S}\\i-j=v}}h_i
   +\sum_{\substack{i\in\bar{Z},\ j\in\bar{S}\\i-j=v}}h_i
   +\sum_{\substack{i,j\in\bar{S}\\i-j=v}}h_j
   +\sum_{\substack{i\in\bar{S},\ j\in\bar{Z}\\i-j=v}}h_j\notag\\
 &=\sum_{\substack{i,j\in\bar{S}\\i-j=v}}(h_i+h_j)
   +\sum_{\substack{i\in\bar{Z},\ j\in\bar{S}\\i-j=v}}h_i
   +\sum_{\substack{i\in\bar{S},\ j\in\bar{Z}\\i-j=v}}h_j,
 \qquad v\in\Gamma.
 \label{eq:component}
\end{align}
In \eqref{eq:component}, the first sum is \(0\) if
\(v\notin\bar{S}-\bar{S}\), and the other two are \(0\) if
\(v\notin(\bar{Z}-\bar{S})\cup(\bar{S}-\bar{Z})\).
Hence
\((A'(\mathbf1_{\bar{S}})h)_v=0\) for
\(v\notin(\bar{S}-\bar{S})\cup(\bar{Z}-\bar{S})\cup(\bar{S}-\bar{Z})\), and
\begin{equation}
 \|A'(\mathbf1_{\bar{S}})h\|_2^2
 =\sum_{v\in(\bar{S}-\bar{S})\cup(\bar{Z}-\bar{S})\cup(\bar{S}-\bar{Z})}
   (A'(\mathbf1_{\bar{S}})h)_v^2.
 \label{eq:norm-support}
\end{equation}
Since \(\bar{S}-\bar{S}=-(\bar{S}-\bar{S})\), condition (ii) gives
\begin{equation}
 (\bar{S}-\bar{S})\cap\bigl((\bar{Z}-\bar{S})\cup(\bar{S}-\bar{Z})\bigr)=\varnothing.
 \label{eq:disjoint}
\end{equation}
Hence \eqref{eq:norm-support} splits as
\[
 \|A'(\mathbf1_{\bar{S}})h\|_2^2
 =\underbrace{\sum_{v\in\bar{S}-\bar{S}}(A'(\mathbf1_{\bar{S}})h)_v^2}_{\text{(I)}}
  +\underbrace{\sum_{v\in(\bar{Z}-\bar{S})\cup(\bar{S}-\bar{Z})}
    (A'(\mathbf1_{\bar{S}})h)_v^2}_{\text{(II)}}.
\]

\paragraph{Step 2.}
We show that \(\text{(I)}\geq2\sum_{i\in\bar{S}}h_i^2\). Let
\(h_{\bar{S}}\in\mathbb R^\Gamma_{\bar{S}}\) equal \(h\) on \(\bar{S}\) and
\(0\) elsewhere. For \(v\in\bar{S}-\bar{S}\), \eqref{eq:component} and
\eqref{eq:disjoint} give
\[
 (A'(\mathbf1_{\bar{S}})h)_v
 =\sum_{\substack{i,j\in\bar{S}\\i-j=v}}(h_i+h_j)
 =(A'(\mathbf1_{\bar{S}})h_{\bar{S}})_v ,
\]
where the last equality above is \eqref{eq:Ah-pairs} with \(\bar{S}\) in place
of \(T\).
Hence, by \eqref{eq:norm-support} with \(\bar{S}\) in place of \(T\),
\(\text{(I)}=\|A'(\mathbf1_{\bar{S}})h_{\bar{S}}\|_2^2\).

Since \(\bar{S}\) is collision-free by (i), \(G(\bar{S})=K_{\bar{S}}\) and
\(\Delta(\overline{G(\bar{S})})=0\). Hence \eqref{eq:quantitative} with
\(\bar{S}\) in place of \(T\), \(x^*=\mathbf1_{\bar{S}}\), \(s=k\) and
\(m=M=1\) gives
\(\sigma_{\min,\bar{S}}(A'(\mathbf1_{\bar{S}}))\geq\sqrt{2(k-2)}\). Since
\(k\geq3\), it follows that
\[
 \text{(I)}
 =\|A'(\mathbf1_{\bar{S}})h_{\bar{S}}\|_2^2
 \geq2(k-2)\|h_{\bar{S}}\|_2^2
 \geq2\sum_{i\in\bar{S}}h_i^2 .
\]

\paragraph{Step 3.}
We show that \(\text{(II)}\geq\sum_{z\in\bar{Z}}h_z^2\). For \(v\in(\bar{Z}-\bar{S})\cup(\bar{S}-\bar{Z})\), \eqref{eq:component}
and \eqref{eq:disjoint} give
\begin{equation}
 (A'(\mathbf1_{\bar{S}})h)_v
 =\sum_{\substack{i\in\bar{Z},\ j\in\bar{S}\\i-j=v}}h_i
 +\sum_{\substack{i\in\bar{S},\ j\in\bar{Z}\\i-j=v}}h_j
 =\sum_{i\in\bar{Z}}\bigl(\mathbf1_{\{i-v\in\bar{S}\}}
   +\mathbf1_{\{i+v\in\bar{S}\}}\bigr)h_i .
 \label{eq:inner-II}
\end{equation}
In \eqref{eq:inner-II}, the coefficient of \(h_i\) vanishes exactly when
\(v\notin(i-\bar{S})\cup(\bar{S}-i)\). By (iii), for each \(z\in\bar{Z}\)
there exists
\begin{equation}
 v_z\in(z-\bar{S})\setminus
 \bigcup_{w\in \bar{Z}\setminus\{z\}}
 \bigl((w-\bar{S})\cup(\bar{S}-w)\bigr).
 \label{eq:vz}
\end{equation}
At \(v=v_z\), the coefficient of \(h_w\) thus vanishes for every
\(w\in\bar{Z}\setminus\{z\}\), and \(z-v_z\in\bar{S}\). Hence
\eqref{eq:inner-II} gives
\[
 (A'(\mathbf1_{\bar{S}})h)_{v_z}
 =\bigl(1+\mathbf1_{\{z+v_z\in\bar{S}\}}\bigr)h_z,
 \qquad z\in\bar{Z}.
\]
By \eqref{eq:vz}, \(v_z\neq v_w\) for distinct \(z,w\in\bar{Z}\), and
\(\{v_z:z\in\bar{Z}\}\subseteq\bar{Z}-\bar{S}\). Therefore
\[
 \text{(II)}
 \geq\sum_{z\in\bar{Z}}(A'(\mathbf1_{\bar{S}})h)_{v_z}^2
 =\sum_{z\in\bar{Z}}\bigl(1+\mathbf1_{\{z+v_z\in\bar{S}\}}\bigr)^2h_z^2
 \geq\sum_{z\in\bar{Z}}h_z^2 .
\]

\paragraph{Step 4.}
Adding the bounds of Steps 2 and 3, we obtain
\[
 \|A'(\mathbf1_{\bar{S}})h\|_2^2
 =\text{(I)}+\text{(II)}
 \geq2\sum_{i\in \bar{S}}h_i^2+\sum_{z\in \bar{Z}}h_z^2
 \geq\|h\|_2^2.
\]
Hence \(\smin(A'(\mathbf1_{\bar{S}}))\geq1>0\), and
Proposition~\ref{prop:generic} shows that \(T\) is nonsingular.
\end{proof}

\subsection{Probability bounds for three conditions}
\label{app:probabilities}

Let \(\mathcal E_1,\mathcal E_2,\mathcal E_3\) be the events that
\((\bar{S},\bar{Z})\) satisfies (i), (ii), (iii) of
Lemma~\ref{lem:construct}, respectively.
The event \(\mathcal E_1\) depends only on \(\bar{S}\).

We now establish the three individual probability bounds and then
combine them in Lemma~\ref{lem:combined-failure}.

\begin{lemma}
\label{lem:collision-free}
Let \(\bar{S}\) be a uniformly random \(k\)-element subset of \(\Gamma\),
where \(k\leq N\). Then
\begin{equation}
 \mathbb P(\mathcal E_1^c)\leq\frac{3k^4}{N}.
 \label{eq:collision-probability}
\end{equation}
\end{lemma}
\begin{proof}
We call \(\mathcal C\subseteq\Gamma\) a collision if
\(\mathcal C=\{a,b\}\cup\{a_1,b_1\}\) for pairs \(\{a,b\}\neq\{a_1,b_1\}\)
with \(a\neq b\), \(a_1\neq b_1\) and \(a-b=\pm(a_1-b_1)\). A collision has
\(4\) or \(3\) points. Condition (i) fails exactly when \(\bar{S}\) contains
a collision, so the union bound gives
\[
 \mathbb P(\mathcal E_1^c)
 \leq\sum_{\mathcal C\text{ collision}}\mathbb P(\mathcal C\subseteq\bar{S})
 =\sum_{\substack{\mathcal C\text{ collision}\\|\mathcal C|=4}}
   \mathbb P(\mathcal C\subseteq\bar{S})
  +\sum_{\substack{\mathcal C\text{ collision}\\|\mathcal C|=3}}
   \mathbb P(\mathcal C\subseteq\bar{S}).
\]
For a fixed collision \(\mathcal C\),
\[
 \mathbb P(\mathcal C\subseteq\bar{S})
 =\frac{\binom{N-|\mathcal C|}{k-|\mathcal C|}}{\binom Nk}
 =\prod_{j=0}^{|\mathcal C|-1}\frac{k-j}{N-j}
 \leq\Bigl(\frac kN\Bigr)^{|\mathcal C|}.
\]
We now count each kind of collision.

If \(|\mathcal C|=4\), the pairs \(\{a,b\}\) and \(\{a_1,b_1\}\) are disjoint, and
\(a-b=\pm(a_1-b_1)\) gives \(b_1=a_1\mp(a-b)\). Thus
\[
 \mathcal C=\{a,\,b,\,a_1,\,a_1\mp(a-b)\},\qquad a,b,a_1\in\Gamma,
\]
so there are at most \(2N^3\) such \(\mathcal C\).

If \(|\mathcal C|=3\), the pairs share exactly one point; relabeling within each
pair if necessary, we may assume \(a=a_1\). Then \(b\neq b_1\), so
\(a-b=-(a-b_1)\), that is, \(b_1=2a-b\). Thus
\[
 \mathcal C=\{a,\,b,\,2a-b\},\qquad a,b\in\Gamma,
\]
so there are at most \(N^2\) such \(\mathcal C\). Therefore
\[
 \mathbb P(\mathcal E_1^c)
 \leq2N^3\Bigl(\frac kN\Bigr)^4+N^2\Bigl(\frac kN\Bigr)^3
 =\frac{2k^4+k^3}{N}
 \leq\frac{3k^4}{N}.
\]
\end{proof}

For the remaining estimates, assume \(k<N\).

\begin{lemma}
\label{lem:cond-ii-prob}
For every realization of \(\bar{S}\),
\begin{equation}
 \mathbb P(\mathcal E_2^c\mid \bar{S})\leq
                  \frac{k^3(t-k)}{N-k}.
 \label{eq:cond-ii-prob}
\end{equation}
\end{lemma}
\begin{proof}
If condition (ii) fails, then \(z-a=b-c\) for some \(z\in\bar{Z}\) and
\(a,b,c\in\bar{S}\), that is, \(a+b-c\in\bar{Z}\). Hence
\[
 \mathcal E_2^c\subseteq\bigcup_{a,b,c\in\bar{S}}\{a+b-c\in\bar{Z}\}.
\]
Given \(\bar{S}\), each point of \(\Gamma\setminus\bar{S}\) lies in
\(\bar{Z}\) with probability \(\frac{t-k}{N-k}\), and no point of
\(\bar{S}\) does. By the union bound over the \(k^3\) triples
\((a,b,c)\),
\[
 \mathbb P(\mathcal E_2^c\mid \bar{S})
 \leq\sum_{a,b,c\in\bar{S}}\mathbb P(a+b-c\in\bar{Z}\mid\bar{S})
 \leq\frac{k^3(t-k)}{N-k}.
\]
\end{proof}

\begin{lemma}
\label{lem:cond-iii-prob}
Let \(r\geq1\) be an integer and \(k=(2+2^d)r\). For every
realization of \(\bar{S}\) satisfying (i),
\begin{equation}
 \mathbb P(\mathcal E_3^c\mid \bar{S})
 \leq N\binom{2k^2}{r}\left(\frac tN\right)^{r+1}.
 \label{eq:cond-iii-prob}
\end{equation}
\end{lemma}
\begin{proof}
Condition (iii) fails exactly when some \(z\in \bar{Z}\) satisfies
\begin{equation}
 z-\bar{S}\subseteq
 \bigcup_{w\in \bar{Z}\setminus\{z\}}
 \bigl((w-\bar{S})\cup(\bar{S}-w)\bigr).
 \label{eq:difference-coverage}
\end{equation}

\paragraph{Step 1.}
For distinct \(z,w\in\Gamma\setminus \bar{S}\), we claim that
\begin{equation}
 \left|(z-\bar{S})\cap\bigl((w-\bar{S})\cup(\bar{S}-w)\bigr)\right|
 \leq2+2^d.
 \label{eq:overlap-count}
\end{equation}
We have
\begin{align*}
 |(z-\bar{S})\cap(w-\bar{S})|
 &=\bigl|\{(a,b):\ a,b\in\bar{S},\ z-a=w-b\}\bigr|
  =\bigl|\{(a,b):\ a,b\in\bar{S},\ a-b=z-w\}\bigr|,\\
 |(z-\bar{S})\cap(\bar{S}-w)|
 &=\bigl|\{(a,b):\ a,b\in\bar{S},\ z-a=b-w\}\bigr|
  =\bigl|\{(a,b):\ a,b\in\bar{S},\ a+b=z+w\}\bigr|,
\end{align*}
and we bound the first count by \(2\) and the second by \(2^d\).

For the first count, suppose that \((a,b)\) is a solution, and let
\((a_1,b_1)\) be any solution. Since \(z\neq w\), we have \(a\neq b\) and
\(a_1\neq b_1\), so condition (i) gives
\[
 a-b=z-w=a_1-b_1
 \quad\Longrightarrow\quad \{a_1,b_1\}=\{a,b\}.
\]
Hence every solution is \((a,b)\) or \((b,a)\), and the first count is at
most \(2\).

For the second count, if there exists a solution \((a,b)\) with
\(a\neq b\), then \((a,b)\) and \((b,a)\) are the only solutions: for another
solution \((a_1,b_1)\) with \(a_1\neq a\), condition (i) gives
\[
 a+b=a_1+b_1=z+w
 \;\Longrightarrow\; a-a_1=b_1-b\neq0
 \;\Longrightarrow\; \{a,a_1\}=\{b_1,b\}
 \;\Longrightarrow\; a_1=b,\ b_1=a,
\]
where the last implication uses \(a\neq b\). Otherwise every solution
has \(a=b\). Then \(a+b=z+w\) becomes \(2a=z+w\), which has at most
\(|\Gamma[2]|\leq2^d\) solutions.
Adding the two bounds proves \eqref{eq:overlap-count}.

\paragraph{Step 2.}
For \(z\in\Gamma\setminus \bar{S}\), let
\[
 D(z)=\left\{w\in\Gamma\setminus(\bar{S}\cup\{z\}):
 (z-\bar{S})\cap\bigl((w-\bar{S})\cup(\bar{S}-w)\bigr)\neq\varnothing
 \right\}.
\]
This set is determined by \(\bar{S}\) and \(z\), before \(\bar{Z}\) is
sampled. We show that
\begin{alignat}{2}
 |D(z)|&\leq2k^2&\qquad&\text{for every }z\in\Gamma\setminus\bar{S},
 \label{eq:D-size}\\
 |D(z)\cap\bar{Z}|&\geq r&&\text{for some }z\in\bar{Z},
 \text{ if (iii) fails}.
 \label{eq:D-hit}
\end{alignat}
For \(w\in D(z)\), there exist \(a,b\in\bar{S}\) such that \(z-a=w-b\) or
\(z-a=b-w\). This gives \(w=z+b-a\) or \(w=a+b-z\). Thus
\[
 D(z)\subseteq\{z+b-a,\ a+b-z:\ a,b\in\bar{S}\},
\]
so \(|D(z)|\leq2k^2\).

For \eqref{eq:D-hit}, suppose that (iii) fails, and let \(z\in\bar{Z}\) satisfy
\eqref{eq:difference-coverage}. Every point of \(z-\bar{S}\) lies in
\((w-\bar{S})\cup(\bar{S}-w)\) for some \(w\in\bar{Z}\setminus\{z\}\), and
such a \(w\) lies in \(D(z)\) by definition. Hence
\eqref{eq:overlap-count} gives
\[
 k=|z-\bar{S}|
 \leq\sum_{w\in D(z)\cap \bar{Z}}
 \left|(z-\bar{S})\cap\bigl((w-\bar{S})\cup(\bar{S}-w)\bigr)\right|
 \leq(2+2^d)\,|D(z)\cap \bar{Z}|,
\]
and \(k=(2+2^d)r\) gives \eqref{eq:D-hit}. 

By \eqref{eq:D-hit}, if (iii)
fails, then \(\bar{Z}\) contains \(\{z\}\cup U\) for some
\(z\in\Gamma\setminus\bar{S}\) and some \(r\)-element \(U\subseteq D(z)\).
Hence
\[
 \mathcal E_3^c\subseteq
 \bigcup_{z\in\Gamma\setminus\bar{S}}\ 
 \bigcup_{\substack{U\subseteq D(z)\\|U|=r}}
 \bigl\{\{z\}\cup U\subseteq\bar{Z}\bigr\}.
\]

\paragraph{Step 3.}
Given \(\bar{S}\), each set \(\{z\}\cup U\) in this union is a fixed set
of \(r+1\) points of \(\Gamma\setminus\bar{S}\). Since \(\bar{Z}\) is a
uniformly random \((t-k)\)-element subset of these \(N-k\) points,
\[
 \mathbb P(\{z\}\cup U\subseteq \bar{Z}\mid \bar{S})
 =\frac{\binom{N-k-r-1}{t-k-r-1}}{\binom{N-k}{t-k}}
 =\prod_{j=0}^r\frac{t-k-j}{N-k-j}
 \leq\Bigl(\frac tN\Bigr)^{r+1}.
\]
There are at most
\(N\) choices of \(z\) and, by \eqref{eq:D-size}, at most
\(\binom{2k^2}{r}\) choices of \(U\) for each \(z\). The union bound
therefore gives
\[
 \mathbb P(\mathcal E_3^c\mid \bar{S})
 \leq\sum_{z\in\Gamma\setminus\bar{S}}\
 \sum_{\substack{U\subseteq D(z)\\|U|=r}}
 \mathbb P(\{z\}\cup U\subseteq\bar{Z}\mid\bar{S})
 \leq N\binom{2k^2}{r}\Bigl(\frac tN\Bigr)^{r+1}.
\]
\end{proof}

\begin{lemma}[Combined failure probability]
\label{lem:combined-failure}
Let \(r\geq1\) be an integer, set \(k=(2+2^d)r<N\), and suppose
\(k\leq t\leq N\). For the random pair \((\bar{S},\bar{Z})\) sampled
as above,
\[
 \mathbb P((\mathcal E_1\cap\mathcal E_2\cap\mathcal E_3)^c)
 \leq\frac{3k^4}{N}
      +\frac{k^3(t-k)}{N-k}
      +N\binom{2k^2}{r}\left(\frac tN\right)^{r+1}.
\]
\end{lemma}
\begin{proof}
We prove the bound by splitting the failure event according to whether
condition (i) holds:
\[
 (\mathcal E_1\cap\mathcal E_2\cap\mathcal E_3)^c
 =\mathcal E_1^c\cup\bigl(\mathcal E_1\cap(\mathcal E_2\cap\mathcal E_3)^c\bigr),
\]
where the union is disjoint. For every \(\bar{S}\in\mathcal E_1\),
Lemmas~\ref{lem:cond-ii-prob} and~\ref{lem:cond-iii-prob} give
\begin{align*}
 \mathbb P((\mathcal E_2\cap\mathcal E_3)^c\mid \bar{S})
 &\leq\mathbb P(\mathcal E_2^c\mid \bar{S})
   +\mathbb P(\mathcal E_3^c\mid \bar{S})\\
 &\leq\frac{k^3(t-k)}{N-k}
       +N\binom{2k^2}{r}\left(\frac tN\right)^{r+1}.
\end{align*}

Lemma~\ref{lem:collision-free} and the conditional bound above give
\begin{align}
 \mathbb P((\mathcal E_1\cap\mathcal E_2\cap\mathcal E_3)^c)
 &=\mathbb P(\mathcal E_1^c)
   +\mathbb P\bigl(\mathcal E_1\cap(\mathcal E_2\cap\mathcal E_3)^c\bigr)\notag\\
 &=\mathbb P(\mathcal E_1^c)
   +\mathbb E\bigl[\mathbf1_{\mathcal E_1}
     \mathbb P((\mathcal E_2\cap\mathcal E_3)^c\mid \bar{S})\bigr]\notag\\
 &\leq\mathbb P(\mathcal E_1^c)
   +\mathbb P(\mathcal E_1)\Bigl(\frac{k^3(t-k)}{N-k}
       +N\binom{2k^2}{r}\Bigl(\frac tN\Bigr)^{r+1}\Bigr)\notag\\
 &\leq\frac{3k^4}{N}
       +\frac{k^3(t-k)}{N-k}
       +N\binom{2k^2}{r}\Bigl(\frac tN\Bigr)^{r+1}.
 \label{eq:generic-failure}
\end{align}
\end{proof}

\subsection{Proof of Theorem~\ref{thm:generic}}
\label{app:generic-proof}
\begin{proof}[Proof of Theorem~\ref{thm:generic}]
By Lemma~\ref{lem:construct},
\begin{equation}
 \mathcal E_1\cap\mathcal E_2\cap\mathcal E_3
 \subseteq\{\bar{S}\cup\bar{Z}\text{ is nonsingular}\}.
 \label{eq:event-inclusion}
\end{equation}
It therefore suffices to bound
\(\mathbb P((\mathcal E_1\cap\mathcal E_2\cap\mathcal E_3)^c)\).

We take \(r,k\) as in the theorem and choose
\(N_*\geq\max\{2k,C^{1/\varepsilon}\}\). For \(N\geq N_*\), this choice
gives \(k<N\) and \(t\leq CN^{1-\varepsilon}\leq N\), so
Lemma~\ref{lem:combined-failure} applies and we obtain
\begin{equation}
 \begin{aligned}
 \mathbb P((\mathcal E_1\cap\mathcal E_2\cap\mathcal E_3)^c)
 &\leq\frac{3k^4}{N}
      +\frac{k^3(t-k)}{N-k}
      +N\binom{2k^2}{r}\left(\frac tN\right)^{r+1}\\
 &\leq3k^4N^{-1}
      +2Ck^3N^{-\varepsilon}
      +C^{r+1}\binom{2k^2}{r}N^{1-\varepsilon(r+1)}\\
 &\leq KN^{-\varepsilon},
 \end{aligned}
 \label{eq:final-failure}
\end{equation}
where \(K=3k^4+2Ck^3+C^{r+1}\binom{2k^2}{r}\) depends
only on \(d,\varepsilon,C\), since \(r,k\) depend only on \(d,\varepsilon\).
Here we use \(N-k\geq\frac N2\), \(t-k\leq t\) and
\(\frac tN\leq CN^{-\varepsilon}\). To obtain the final estimate, we also
use \(N^{-1}\leq N^{-\varepsilon}\) and
\(N^{1-\varepsilon(r+1)}\leq N^{-\varepsilon}\), which hold since
\(0<\varepsilon\leq1\) and \(\varepsilon r\geq1\).

Finally, \eqref{eq:same-law}, \eqref{eq:event-inclusion} and
\eqref{eq:final-failure} give
\[
 \mathbb P(T\text{ is nonsingular})
 =\mathbb P(\bar{S}\cup\bar{Z}\text{ is nonsingular})
 \geq\mathbb P(\mathcal E_1\cap\mathcal E_2\cap\mathcal E_3)
 \geq1-KN^{-\varepsilon}.
\]
If \(T\) is nonsingular, then \(\smin(A'(x))>0\) for almost every
\(x\in\RT\), and Lemma~\ref{lem:local} gives local Lipschitz stability
at each of these points.
\end{proof}

\section{Proof of Lemma~\ref{lem:random-degree}}
\label{app:degree}

Lemma~\ref{lem:random-degree} fails only if one of the events
\[
 \bigl\{\bigl||T|-t\bigr|\geq\theta t\bigr\},
 \qquad
 \bigl\{\Delta(\overline{G(T)})>\theta t\bigr\}
\]
occurs. For the first event, \(|T|\sim\operatorname{Bin}(N,\frac tN)\) has
mean \(t\). Since \(\theta<\frac12\), combining the
upper and lower Chernoff tails
in~\cite[Theorem~1.1, inequality~(1.7)]{DubhashiPanconesi2009} gives
\begin{equation}
 \mathbb P\bigl(\bigl||T|-t\bigr|\geq\theta t\bigr)
 \leq2e^{-\theta^2t/3}.
 \label{eq:size-tail}
\end{equation}
For the second event, we prove that, for sufficiently large \(N\),
\begin{equation}
 \mathbb P\left(\Delta(\overline{G(T)})>\theta t\right)
 \leq C_1t e^{-c_1\sqrt t\log(N/t)}+C_2t e^{-c_2t},
 \label{eq:max-tail}
\end{equation}
with \(C_1,c_1,C_2,c_2>0\) depending only on \(d,c,\theta\). A union
bound over the two events then gives
\begin{align*}
 \mathbb P\left(\text{Lemma~\ref{lem:random-degree} fails}\right)
 &\leq\mathbb P\bigl(\bigl||T|-t\bigr|\geq\theta t\bigr)
   +\mathbb P\left(\Delta(\overline{G(T)})>\theta t\right)\\
 &\leq2e^{-\theta^2t/3}
   +C_1t e^{-c_1\sqrt t\log(N/t)}+C_2t e^{-c_2t}
  =\varepsilon_\theta(t,N).
\end{align*}
It therefore remains to prove \eqref{eq:max-tail}, which occupies the
remainder of this appendix. In
Appendix~\ref{app:degree-bound}, we bound the degree of each
\(a\in T\) in \(\overline{G(T)}\) by a sum \(W_a\) of products of
indicators. We then split each \(W_a\) into two parts and use a union
bound over \(a\in\Gamma\) to reduce \eqref{eq:max-tail} to tail bounds for
these parts conditional on \(a\in T\). We prove these tail
bounds in Appendices~\ref{app:degree-nd} and~\ref{app:degree-deg}, and
complete the proof of \eqref{eq:max-tail} in
Appendix~\ref{app:degree-completion}.

\subsection{A bound for \texorpdfstring{\(\Delta(\overline{G(T)})\)}{the maximum degree of the complement of G(T)}}
\label{app:degree-bound}
A union bound over the vertices in \(T\) gives
\begin{equation}
 \mathbb P\left(\Delta(\overline{G(T)})>\theta t\right)
 \leq\sum_{a\in\Gamma}
   \mathbb P\left(a\in T,\ d_{\overline{G(T)}}(a)>\theta t\right).
 \label{eq:degree-union}
\end{equation}
 We first bound \(d_{\overline{G(T)}}(a)\) for a fixed \(a\in T\) by a
sum of products of the indicators
\(\xi_x=\xi_x(T)=\mathbf1_{\{x\in T\}}\), \(x\in\Gamma\). Here each \(\xi_x\) is
a function of the random set \(T\), and these indicators are independent
Bernoulli variables of mean \(p=\frac tN\).

We fix \(a\in T\). If \(b\in T\setminus\{a\}\) and
\(\{a,b\}\notin E(T)\), then by the definition of \(E(T)\) in
Subsection~\ref{subsec:isolated-graph}, there are distinct
\(a_1,b_1\in T\) with \(\{a_1,b_1\}\neq\{a,b\}\) and \(a_1-b_1=\pm(a-b)\).
Ordering them so that \(a_1-b_1=a-b\) and setting \(h=a_1-a\), we get
\(\{a_1,b_1\}=\{a+h,b+h\}\neq\{a,b\}\) with \(b,a+h,b+h\in T\). Hence
\[
 \mathbf1_{\{b\in T,\ \{a,b\}\notin E(T)\}}
 \leq\sum_{\substack{h\in\Gamma\\\{a+h,b+h\}\neq\{a,b\}}}
       \xi_b\xi_{a+h}\xi_{b+h}.
\]
Summing over \(b\neq a\), we obtain
\begin{equation}
 d_{\overline{G(T)}}(a)\leq W_a,
 \qquad
 W_a=\sum_{b\neq a}
       \sum_{\substack{h\in\Gamma\\\{a+h,b+h\}\neq\{a,b\}}}
       \xi_b\xi_{a+h}\xi_{b+h}.
 \label{eq:Wa}
\end{equation}

Since \(\{a+h,b+h\}\neq\{a,b\}\), the two pairs share at most one point.
We therefore split \(W_a=W_a^{\mathrm{nd}}+W_a^{\mathrm{deg}}\) according
to whether they are disjoint or share exactly one point:
\begin{equation}
 W_a^{\mathrm{nd}}
 =\sum_{b\neq a}
   \sum_{\substack{h\in\Gamma\\\{a+h,b+h\}\cap\{a,b\}=\varnothing}}
   \xi_b\xi_{a+h}\xi_{b+h},
 \qquad
 W_a^{\mathrm{deg}}
 =\sum_{b\neq a}
   \sum_{\substack{h\in\Gamma\\|\{a+h,b+h\}\cap\{a,b\}|=1}}
   \xi_b\xi_{a+h}\xi_{b+h}.
 \label{eq:Wa-split}
\end{equation}

We fix \(\eta=\frac{\theta-6c^2}{8}\), which is positive since
\(c^2<\frac\theta6\). By \eqref{eq:degree-union} and \eqref{eq:Wa}, we
have
\begin{equation}
 \begin{aligned}
 \mathbb P\left(\Delta(\overline{G(T)})>\theta t\right)
 &\leq\sum_{a\in\Gamma}\mathbb P\left(a\in T,\ W_a\geq\theta t\right)
  =p\sum_{a\in\Gamma}\mathbb P\left(W_a\geq\theta t\,\middle|\,a\in T\right)\\
 &\leq p\sum_{a\in\Gamma}\Bigl[
   \mathbb P\bigl(W_a^{\mathrm{nd}}\geq6(c^2+\eta)t\mid a\in T\bigr)
   +\mathbb P\bigl(W_a^{\mathrm{deg}}\geq2\eta t\mid a\in T\bigr)\Bigr].
 \end{aligned}
 \label{eq:max-reduction}
\end{equation}
For the final estimate, we use \(6(c^2+\eta)+2\eta=\theta\): if
\(W_a^{\mathrm{nd}}<6(c^2+\eta)t\) and \(W_a^{\mathrm{deg}}<2\eta t\),
then \(W_a<\theta t\). To obtain \eqref{eq:max-tail}, it remains to
estimate the two probabilities on the right, which we do in
Appendices~\ref{app:degree-nd} and~\ref{app:degree-deg}.

\subsection{Estimation of \texorpdfstring{\(W_a^{\mathrm{nd}}\)}{W\_a\textasciicircum nd}}
\label{app:degree-nd}
We fix \(a\in\Gamma\) and condition on \(a\in T\), so that \(\xi_a=1\).
We prove that, for sufficiently large \(N\),
\begin{equation}
 \mathbb P\bigl(W_a^{\mathrm{nd}}\geq6(c^2+\eta)t\mid a\in T\bigr)
 \leq C_1e^{-c_1\sqrt t\log(N/t)},
 \label{eq:nd-tail}
\end{equation}
where \(C_1,c_1>0\) depend only on \(d,c,\theta\).
We prove it using the following upper-tail bound, which is the
\(3\)-uniform case
of~\cite[Theorem~6, condition~(b), inequality~(9)]{Warnke2017UpperTails}.\footnote{The
parameters \(a,t,c\) of~\cite{Warnke2017UpperTails} are renamed
\(a_0,\tau,c_0\) to avoid conflicts with our notation.}

\begin{lemma}[Warnke's upper-tail bound]
\label{lem:warnke}
Fix \(a_0>0\), \(D\geq1\), and \(\gamma\in(0,1)\). Let \(\mathcal H\) be a
hypergraph, consisting of a finite vertex set \(V(\mathcal H)\) and a set
\(E(\mathcal H)\) of subsets of \(V(\mathcal H)\), called hyperedges.
Assume that \(\mathcal H\) is \(3\)-uniform, that is, \(|f|=3\) for every
\(f\in E(\mathcal H)\). Write
\(v(\mathcal H)=|V(\mathcal H)|\) for the number of vertices,
\(e(\mathcal H)=|E(\mathcal H)|\) for the number of hyperedges, and
\[
 \Delta_2(\mathcal H)
 =\max_{\substack{S\subseteq V(\mathcal H)\\|S|=2}}
   |\{f\in E(\mathcal H):S\subseteq f\}|
\]
for the maximum pair-codegree. Suppose that
\[
 v(\mathcal H)\leq Dn,\qquad e(\mathcal H)\geq a_0n^2,\qquad
 \Delta_2(\mathcal H)\leq D.
\]
Let \(R\subseteq V(\mathcal H)\) be a random set that contains each vertex
independently with probability \(p\in(0,1]\). Let
\(\xi_x(R)=\mathbf1_{\{x\in R\}}\) for
\(x\in V(\mathcal H)\), and set
\[
 \begin{aligned}
 X&=\sum_{f\in E(\mathcal H)}\prod_{x\in f}\xi_x(R)
   =|\{f\in E(\mathcal H):f\subseteq R\}|,
 &\mu&=\mathbb EX=e(\mathcal H)p^3,\\
 \Lambda&=\mu(1+np^2),
 &\varphi(u)&=(1+u)\log(1+u)-u.
 \end{aligned}
\]
There exist constants \(n_0=n_0(a_0,D)>0\) and
\(c_0=c_0(a_0,D,\gamma)>0\) such that, for every \(n\geq n_0\) and
\(\tau>0\) satisfying
\[
 \tau\geq\gamma\min\{(\Var X)^{2/3},\mu^{2/3}\}(\log n)^{4/3},
\]
we have
\begin{equation}
 \mathbb P(X\geq\mu+\tau)
 \leq(1+n^{-1})\exp\left[-c_0\min\left\{
          \frac{\varphi(\tau/\mu)\mu^2}{\Lambda},
          \sqrt\tau\log(e/p)\right\}\right].
 \label{eq:warnke}
\end{equation}
\end{lemma}

\paragraph{Step 1.}
We first replace the sum over ordered triples \((b,b+h,a+h)\) in
\(W_a^{\mathrm{nd}}\) by a sum over the unordered hyperedges of a hypergraph.
In each term of \(W_a^{\mathrm{nd}}\), 
the pairs \(\{a,b\}\) and \(\{a+h,b+h\}\) are disjoint, so the points
\(b,a+h,b+h\) are distinct and differ from \(a\). We let \(\mathcal H_a\)
be the \(3\)-uniform hypergraph on \(\Gamma\) with hyperedge set
\[
 E(\mathcal H_a)=\bigl\{\{b,b+h,a+h\}:\ b\neq a,\ h\neq0,\ b\neq a+h,\
 b+h\neq a\bigr\},
\]
so that \(a\) is an isolated vertex of \(\mathcal H_a\), and define
\[
 X_a=\sum_{f\in E(\mathcal H_a)}\prod_{x\in f}\xi_x
 =\sum_{\{b,b+h,a+h\}\in E(\mathcal H_a)}\xi_b\xi_{a+h}\xi_{b+h}.
\]
Comparing with \eqref{eq:Wa-split}, \(W_a^{\mathrm{nd}}\) sums the same
terms over the ordered triples \((b,b+h,a+h)\), and each unordered hyperedge of
\(\mathcal H_a\) has at least one and at most \(3!=6\) such orderings. Hence
\begin{equation}
 X_a\leq W_a^{\mathrm{nd}}\leq6X_a,
 \qquad\text{so}\qquad
 \{W_a^{\mathrm{nd}}\geq6(c^2+\eta)t\}
 \subseteq\{X_a\geq(c^2+\eta)t\}.
 \label{eq:nd-inclusion}
\end{equation}

\paragraph{Step 2.}
We show that, for \(N\geq8\), \(\mathcal H_a\) satisfies the hypergraph
conditions of Lemma~\ref{lem:warnke} with \(n=N\), \(a_0=\frac1{12}\), and
\(D=3\), that is,
\begin{equation}
 v(\mathcal H_a)\leq N,
 \qquad
 \frac{N^2}{12}\leq e(\mathcal H_a)\leq N^2,
 \qquad
 \Delta_2(\mathcal H_a)\leq3.
 \label{eq:hypergraph-size}
\end{equation}
\begin{enumerate}[label=(\roman*)]
\item The first bound holds since \(\mathcal H_a\) has vertex set
\(\Gamma\).
\item For the upper bound, each unordered hyperedge \(\{b,b+h,a+h\}\) of
\(\mathcal H_a\) comes from at least one ordered pair \((b,h)\), so
\[
 e(\mathcal H_a)\leq\bigl|\{(b,h):b\neq a,\ h\neq0\}\bigr|=(N-1)^2\leq N^2.
\]
For the lower bound, by Step~1, each unordered hyperedge of \(\mathcal H_a\)
comes from at most six ordered pairs \((b,h)\), so
\[
 e(\mathcal H_a)
 \geq\frac16\bigl|\{(b,h):b\neq a,\ h\neq0,\ b\neq a+h,\ b+h\neq a\}\bigr|
 \geq\frac{(N-1)^2-2(N-1)}6
 \geq\frac{N^2}{12},
\]
where the conditions \(b\neq a\) and \(h\neq0\) give exactly
\((N-1)^2\) ordered pairs \((b,h)\), and the constraints \(b\neq a+h\) and
\(b+h\neq a\) each exclude at most \(N-1\) of them. The final estimate
uses \(N\geq8\).
\item By the definition in Lemma~\ref{lem:warnke},
\[
 \Delta_2(\mathcal H_a)
 =\max_{\substack{x,y\in\Gamma\\x\neq y}}
   \bigl|\{z:\{x,y,z\}\in E(\mathcal H_a)\}\bigr|
 \leq\max_{\substack{x,y\in\Gamma\\x\neq y}}
   \bigl|\{x+y-a,\ a+x-y,\ a+y-x\}\bigr|
 \leq3,
\]
since each hyperedge \(\{b,b+h,a+h\}\) satisfies \(a+(b+h)=b+(a+h)\), so
\(\{x,y,z\}\in E(\mathcal H_a)\) implies \(a+z=x+y\), \(a+x=y+z\), or
\(a+y=x+z\).
\end{enumerate}

\paragraph{Step 3.}
We apply Lemma~\ref{lem:warnke} to
\(X_a=\sum_{f\in E(\mathcal H_a)}\prod_{x\in f}\xi_x(T)\) with \(R=T\),
\(\gamma=\frac12\), and \(\tau_a\) defined by
\[
 (c^2+\eta)t=\mu_a+\tau_a,
 \qquad
 \mu_a=\mathbb EX_a=e(\mathcal H_a)p^3.
\]
Since \(p=\frac tN\) and \(t\leq c\sqrt N\), \eqref{eq:hypergraph-size}
gives
\begin{equation}
 \mu_a\leq N^2p^3=\frac{t^3}{N}\leq c^2t,
 \qquad
 \tau_a=(c^2+\eta)t-\mu_a\geq\eta t.
 \label{eq:mu-tau}
\end{equation}
Since \(\tau_a\geq\eta t\), the condition on \(\tau\) in
Lemma~\ref{lem:warnke} holds for sufficiently large \(N\), because
\((\log N)^4=o(t)\) gives
\[
 \frac12\min\{(\Var X_a)^{2/3},\mu_a^{2/3}\}(\log N)^{4/3}
 \leq\frac12c^{4/3}t^{2/3}(\log N)^{4/3}
 =\frac12c^{4/3}t\left(\frac{(\log N)^4}{t}\right)^{1/3}
 =o(t).
\]
Hence \eqref{eq:nd-inclusion} and \eqref{eq:warnke} give, for sufficiently large \(N\),
\begin{equation}
 \begin{aligned}
 \mathbb P(W_a^{\mathrm{nd}}\geq6(c^2+\eta)t\mid a\in T)
 &\leq\mathbb P(X_a\geq\mu_a+\tau_a\mid a\in T)
  =\mathbb P(X_a\geq\mu_a+\tau_a)\\
 &\leq(1+N^{-1})\exp\left[-c_0\min\left\{
       \frac{\varphi(\tau_a/\mu_a)\mu_a^2}{\Lambda},
       \sqrt{\tau_a}\log(e/p)\right\}\right].
 \end{aligned}
 \label{eq:nd-warnke}
\end{equation}

\paragraph{Step 4.}
We simplify the right-hand side of \eqref{eq:nd-warnke}. On the one hand,
since \(\Lambda=\mu_a(1+Np^2)\) and \(\varphi(u)=(1+u)\log(1+u)-u\), we
have
\[
 \frac{\varphi(\tau_a/\mu_a)\mu_a^2}{\Lambda}
 =\frac{\tau_a}{1+Np^2}\,\frac{\varphi(\tau_a/\mu_a)}{\tau_a/\mu_a}
 \geq\frac{c^2\varphi(\eta/c^2)}{1+c^2}\,t.
\]
Here we use \(\frac{\tau_a}{\mu_a}\geq\frac{\eta}{c^2}\) from
\eqref{eq:mu-tau}, and that \(\frac{\varphi(u)}u\) is increasing for
\(u>0\).

On the other hand, since \(\log(e/p)=1+\log(N/t)\), we have
\[
 \sqrt{\tau_a}\log(e/p)
 =\sqrt{\tau_a}\bigl(1+\log(N/t)\bigr)
 \geq\sqrt{\eta t}\log(N/t).
\]
Here we use \(\tau_a\geq\eta t\) from \eqref{eq:mu-tau}. As
\(\sqrt t\log(N/t)\leq\sqrt t\log N=o(t)\), for sufficiently large \(N\),
\[
 \min\left\{\frac{\varphi(\tau_a/\mu_a)\mu_a^2}{\Lambda},
            \sqrt{\tau_a}\log(e/p)\right\}
 \geq\min\left\{\frac{c^2\varphi(\eta/c^2)}{1+c^2}\,t,\
            \sqrt{\eta t}\log(N/t)\right\}
 =\sqrt{\eta t}\log(N/t).
\]
Substituting this into \eqref{eq:nd-warnke} and using \(1+N^{-1}\leq2\), we
obtain
\[
 \mathbb P(W_a^{\mathrm{nd}}\geq6(c^2+\eta)t\mid a\in T)
 \leq2e^{-c_0\sqrt\eta\sqrt t\log(N/t)},
\]
which proves \eqref{eq:nd-tail} with \(C_1=2\) and \(c_1=c_0\sqrt\eta\).

\subsection{Estimation of \texorpdfstring{\(W_a^{\mathrm{deg}}\)}{W\_a\textasciicircum deg}}
\label{app:degree-deg}
We continue to fix \(a\) and condition on \(a\in T\). We prove that,
for sufficiently large \(N\),
\begin{equation}
 \mathbb P(W_a^{\mathrm{deg}}\geq2\eta t\mid a\in T)
 \leq C_2e^{-c_2t},
 \label{eq:deg-tail}
\end{equation}
with \(C_2,c_2>0\) depending only on \(d,c,\theta\).

\paragraph{Step 1.}
We first show that the double sum in \eqref{eq:Wa-split} reduces to a sum
over \(b\) alone:
\begin{equation}
 W_a^{\mathrm{deg}}
 =\sum_{b\neq a}
   \sum_{\substack{h\in\Gamma\\|\{a+h,b+h\}\cap\{a,b\}|=1}}
   \xi_b\xi_{a+h}\xi_{b+h}
 =\sum_{\substack{b\in\Gamma\\2(b-a)\neq0}}
   \bigl(\xi_b\xi_{2b-a}+\xi_b\xi_{2a-b}\bigr).
 \label{eq:deg-b-sum}
\end{equation}
For \(b\neq a\) and \(h\in\Gamma\), the pairs \(\{a,b\}\) and
\(\{a+h,b+h\}\) share exactly one point if and only if one of the
following cases holds:
\[
 \{a+h,b+h\}=
 \begin{cases}
  \{b,2b-a\}, & h=b-a,\ 2b-a\neq a,\\
  \{2a-b,a\}, & h=a-b,\ 2a-b\neq b.
 \end{cases}
\]
Both \(2b-a\neq a\) and \(2a-b\neq b\) mean \(2(b-a)\neq0\). Hence, for
\(b\neq a\), \(|\{a+h,b+h\}\cap\{a,b\}|=1\) if and only if \(h\in\{b-a,a-b\}\) and
\(2(b-a)\neq0\). Therefore,
\[
 W_a^{\mathrm{deg}}
 =\sum_{\substack{b\in\Gamma\\2(b-a)\neq0}}
   \bigl(\xi_b\xi_b\xi_{2b-a}+\xi_b\xi_{2a-b}\xi_a\bigr)
 =\sum_{\substack{b\in\Gamma\\2(b-a)\neq0}}
   \bigl(\xi_b\xi_{2b-a}+\xi_b\xi_{2a-b}\bigr),
\]
where we used \(\xi_b^2=\xi_b\) and \(\xi_a=1\). 

\paragraph{Step 2.}
We next replace the sums over ordered pairs in \(W_a^{\mathrm{deg}}\) by a
sum over the unordered edges of a graph. We let \(G_a\) be the simple graph
on \(\Gamma\setminus\{a\}\) with edge set
\[
 E(G_a)=\bigl\{\{b,2b-a\},\ \{b,2a-b\}:\ b\in\Gamma,\ 2(b-a)\neq0\bigr\}.
\]
The terms of \eqref{eq:deg-b-sum} are indexed by the ordered pairs
\((b,2b-a)\) and \((b,2a-b)\), and \(E(G_a)\) consists of the
corresponding unordered edges. An edge \(\{u,v\}\) corresponds to the
terms with
\[
 (b,2b-a)=(u,v),\qquad
 (b,2b-a)=(v,u),\qquad
 (b,2a-b)=(u,v),\qquad
 (b,2a-b)=(v,u).
\]
These equations force \(b=u\), \(b=v\), \(b=u\) and \(b=v\), respectively,
so each edge corresponds to at most four terms. Therefore
\begin{equation}
 W_a^{\mathrm{deg}}
 \leq4\sum_{\{u,v\}\in E(G_a)}\xi_u\xi_v.
 \label{eq:deg-pairs}
\end{equation}

\paragraph{Step 3.}
We partition \(E(G_a)\) into matchings\footnote{A matching is a
set of edges, no two of which share a vertex.} and reduce
\eqref{eq:deg-tail} to tail bounds for the matching sums.
By Vizing's theorem~\cite{Vizing1964ChromaticClass}, the edge set of a
simple graph
\(G\) can be partitioned into at most \(\Delta(G)+1\) matchings. We next
bound \(\Delta(G_a)\).

Let \(u\in\Gamma\setminus\{a\}\). By the definition of \(E(G_a)\), every
edge containing \(u\) is \(\{b,2b-a\}\) or \(\{b,2a-b\}\) for some
\(b\in\Gamma\), with \(u\) as an endpoint. Depending on its form and on
which endpoint equals \(u\), the edge is
\[
 \begin{cases}
  \{u,2u-a\}, & \{b,2b-a\}\ \text{with}\ b=u,\\
  \{b,u\}, & \{b,2b-a\}\ \text{with}\ 2b-a=u,\\
  \{u,2a-u\}, & \{b,2a-b\}\ \text{with}\ b=u\ \text{or}\ 2a-b=u.
 \end{cases}
\]
In the second case, the equation \(2b-a=u\) has at most
\(|\Gamma[2]|\leq2^d\) solutions \(b\). Hence \(d_{G_a}(u)\leq2+2^d\), so \(\Delta(G_a)\leq2+2^d\), and by
Vizing's theorem \(E(G_a)\) can be partitioned into \(2^d+3\) matchings
\(E_1,\dots,E_{2^d+3}\). Then
\eqref{eq:deg-pairs} gives
\begin{equation}
 \{W_a^{\mathrm{deg}}\geq2\eta t\}
 \subseteq
 \left\{\sum_{\ell=1}^{2^d+3}
       \sum_{\{u,v\}\in E_\ell}\xi_u\xi_v\geq\frac{\eta t}{2}\right\}
 \subseteq\bigcup_{\ell=1}^{2^d+3}
 \left\{\sum_{\{u,v\}\in E_\ell}\xi_u\xi_v
                    \geq\frac{\eta t}{2(2^d+3)}\right\}.
 \label{eq:deg-inclusion}
\end{equation}

\paragraph{Step 4.}
Finally, we bound each matching sum by the Chernoff bound and combine the
bounds by a union bound. We first prove that, for
sufficiently large \(N\) and every \(\ell\),
\begin{equation}
 \mathbb P\left(\sum_{\{u,v\}\in E_\ell}\xi_u\xi_v
   \geq\frac{\eta t}{2(2^d+3)}\,\middle|\,a\in T\right)
 \leq\exp\left(-\frac{\eta t}{2(2^d+3)}\right).
 \label{eq:matching-tail}
\end{equation}
By the Chernoff bound~\cite[Sections~1.3--1.4]{DubhashiPanconesi2009}, if
\(X\) is a sum of independent Bernoulli variables, \(\mu=\mathbb EX\), and
\(\kappa\geq e^2\mu\) is a real threshold, then
\[
 \mathbb P(X\geq\kappa)\leq\left(\frac{e\mu}{\kappa}\right)^{\kappa}
 \leq e^{-\kappa}.
\]
We take
\[
 X=\sum_{\{u,v\}\in E_\ell}\xi_u\xi_v,
 \qquad
 \kappa=\frac{\eta t}{2(2^d+3)}.
\]
Within \(E_\ell\), the products \(\xi_u\xi_v\) depend on disjoint sets of
indicators outside \(\{a\}\). So, conditional on \(a\in T\), \(X\) is a sum
of independent Bernoulli variables, and its mean is
\(\mu=\mathbb E[X\mid a\in T]=|E_\ell|p^2\leq Np^2=\frac{t^2}{N}\leq c^2\).
Thus \(\kappa\geq e^2c^2\geq e^2\mu\) for sufficiently large \(N\), and the
Chernoff bound, under the conditional probability given \(a\in T\), gives
\eqref{eq:matching-tail}.

Combining \eqref{eq:deg-inclusion} and \eqref{eq:matching-tail} by a
union bound, we have
\begin{align*}
 \mathbb P(W_a^{\mathrm{deg}}\geq2\eta t\mid a\in T)
 &\leq\sum_{\ell=1}^{2^d+3}
   \mathbb P\left(\sum_{\{u,v\}\in E_\ell}\xi_u\xi_v
   \geq\frac{\eta t}{2(2^d+3)}\,\middle|\,a\in T\right)\\
 &\leq(2^d+3)\exp\left(-\frac{\eta t}{2(2^d+3)}\right),
\end{align*}
which proves \eqref{eq:deg-tail} with \(C_2=2^d+3\) and
\(c_2=\frac{\eta}{2(2^d+3)}\).

\subsection{Completion of the proof of \texorpdfstring{\eqref{eq:max-tail}}{the bound on the maximum degree}}
\label{app:degree-completion}
Substituting \eqref{eq:nd-tail} and \eqref{eq:deg-tail} into
\eqref{eq:max-reduction} and using \(Np=t\), we obtain
\[
 \mathbb P\left(\Delta(\overline{G(T)})>\theta t\right)
 \leq Np\left(C_1e^{-c_1\sqrt t\log(N/t)}+C_2e^{-c_2t}\right)
 =C_1t e^{-c_1\sqrt t\log(N/t)}+C_2t e^{-c_2t},
\]
which is \eqref{eq:max-tail}.

\bibliographystyle{plain}
\bibliography{revised-v5-ref}

\begin{thebibliography}{10}

\bibitem{BalanZou2016}
Radu Balan and Dongmian Zou.
\newblock On {Lipschitz} analysis and {Lipschitz} synthesis for the phase retrieval problem.
\newblock {\em Linear Algebra and its Applications}, 496:152--181, 2016.

\bibitem{Bandeira2014SavingPhase}
Afonso~S. Bandeira, Jameson Cahill, Dustin~G. Mixon, and Aaron~A. Nelson.
\newblock Saving phase: Injectivity and stability for phase retrieval.
\newblock {\em Applied and Computational Harmonic Analysis}, 37(1):106--125, 2014.

\bibitem{barnett2022geometry}
Alexander~H Barnett, Charles~L Epstein, Leslie Greengard, and Jeremy Magland.
\newblock {\em Geometry of the phase retrieval problem: graveyard of algorithms}.
\newblock Cambridge University Press, 2022.

\bibitem{barnett2020geometry}
Alexander~H Barnett, Charles~L Epstein, Leslie~F Greengard, and Jeremy~F Magland.
\newblock Geometry of the phase retrieval problem.
\newblock {\em Inverse Problems}, 36(9):094003, 2020.

\bibitem{BendoryEdidin2020}
Tamir Bendory and Dan Edidin.
\newblock Toward a mathematical theory of the crystallographic phase retrieval problem.
\newblock {\em SIAM Journal on Mathematics of Data Science}, 2(3):809--839, 2020.

\bibitem{Candes2013}
Emmanuel~J. Cand{\`e}s, Thomas Strohmer, and Vladislav Voroninski.
\newblock Phaselift: Exact and stable signal recovery from magnitude measurements via convex programming.
\newblock {\em Communications on Pure and Applied Mathematics}, 66(8):1241--1274, 2013.

\bibitem{DaintyFienup1987}
J.~C. Dainty and J.~R. Fienup.
\newblock Phase retrieval and image reconstruction for astronomy.
\newblock In H.~Stark, editor, {\em Image Recovery: Theory and Application}, chapter~7, pages 231--275. Academic Press, 1987.

\bibitem{DubhashiPanconesi2009}
Devdatt~P. Dubhashi and Alessandro Panconesi.
\newblock {\em Concentration of Measure for the Analysis of Randomized Algorithms}.
\newblock Cambridge University Press, 2009.

\bibitem{EdidinSuresh2026}
Dan Edidin and Arun Suresh.
\newblock The generic crystallographic phase retrieval problem.
\newblock {\em Applied and Computational Harmonic Analysis}, 84:101888, 2026.

\bibitem{EldarMendelson2014}
Yonina~C. Eldar and Shahar Mendelson.
\newblock Phase retrieval: Stability and recovery guarantees.
\newblock {\em Applied and Computational Harmonic Analysis}, 36(3):473--494, 2014.

\bibitem{elser2003phase}
Veit Elser.
\newblock Phase retrieval by iterated projections.
\newblock {\em Journal of the Optical Society of America A}, 20(1):40--55, 2003.

\bibitem{ElserLanBendory2018Benchmark}
Veit Elser, Ti-Yen Lan, and Tamir Bendory.
\newblock Benchmark problems for phase retrieval.
\newblock {\em SIAM Journal on Imaging Sciences}, 11(4):2429--2455, 2018.

\bibitem{fienup1986stagnation}
J.~R. Fienup and C.~C. Wackerman.
\newblock Phase-retrieval stagnation problems and solutions.
\newblock {\em Journal of the Optical Society of America A}, 3(11):1897--1907, 1986.

\bibitem{Fienup1982}
James~R. Fienup.
\newblock Phase retrieval algorithms: a comparison.
\newblock {\em Applied Optics}, 21(15):2758--2769, 1982.

\bibitem{Gerchberg1972}
R.~W. Gerchberg and W.~O. Saxton.
\newblock A practical algorithm for the determination of phase from image and diffraction plane pictures.
\newblock {\em Optik}, 35(2):237--246, 1972.

\bibitem{ghosh2023sparse}
Subhroshekhar Ghosh and Philippe Rigollet.
\newblock Sparse multi-reference alignment: Phase retrieval, uniform uncertainty principles and the beltway problem.
\newblock {\em Foundations of Computational Mathematics}, 23(5):1851--1898, 2023.

\bibitem{Jaganathan2015Overview}
Kishore Jaganathan, Yonina~C. Eldar, and Babak Hassibi.
\newblock Phase retrieval: An overview of recent developments.
\newblock In Adrian Stern, editor, {\em Optical Compressive Imaging}, pages 263--296. CRC Press, Boca Raton, FL, 2016.

\bibitem{jaganathan2017sparse}
Kishore Jaganathan, Samet Oymak, and Babak Hassibi.
\newblock Sparse phase retrieval: Uniqueness guarantees and recovery algorithms.
\newblock {\em IEEE Transactions on Signal Processing}, 65(9):2402--2410, 2017.

\bibitem{Luke2005}
D.~Russell Luke.
\newblock Relaxed averaged alternating reflections for diffraction imaging.
\newblock {\em Inverse Problems}, 21(1):37--50, 2005.

\bibitem{Luke2017WhatsNew}
D.~Russell Luke.
\newblock Phase retrieval, what's new?
\newblock {\em SIAG/OPT Views and News}, 25(1):1--6, 2017.

\bibitem{Marchesini2007}
Stefano Marchesini.
\newblock A unified evaluation of iterative projection algorithms for phase retrieval.
\newblock {\em Review of Scientific Instruments}, 78:011301, 2007.

\bibitem{Miao1999}
Jianwei Miao, Pambos Charalambous, Janos Kirz, and David Sayre.
\newblock Extending the methodology of {X-ray} crystallography to allow imaging of micrometre-sized non-crystalline specimens.
\newblock {\em Nature}, 400:342--344, 1999.

\bibitem{Millane1990}
R.~P. Millane.
\newblock Phase retrieval in crystallography and optics.
\newblock {\em Journal of the Optical Society of America A}, 7(3):394--411, 1990.

\bibitem{Ranieri2013}
Juri Ranieri, Amina Chebira, Yue~M. Lu, and Martin Vetterli.
\newblock Phase retrieval for sparse signals: Uniqueness conditions.
\newblock arXiv:1308.3058, 2013.

\bibitem{Shechtman2014GESPAR}
Yoav Shechtman, Amir Beck, and Yonina~C. Eldar.
\newblock Gespar: Efficient phase retrieval of sparse signals.
\newblock {\em IEEE Transactions on Signal Processing}, 62(4):928--938, 2014.

\bibitem{Shechtman2015Overview}
Yoav Shechtman, Yonina~C. Eldar, Oren Cohen, Henry~N. Chapman, Jianwei Miao, and Mordechai Segev.
\newblock Phase retrieval with application to optical imaging: A contemporary overview.
\newblock {\em IEEE Signal Processing Magazine}, 32(3):87--109, 2015.

\bibitem{Vizing1964ChromaticClass}
Vadim~G. Vizing.
\newblock On an estimate of the chromatic class of a {$p$}-graph.
\newblock {\em Diskretnyi Analiz}, 3:25--30, 1964.
\newblock In Russian.

\bibitem{Waldspurger2015}
Ir{\`e}ne Waldspurger, Alexandre d'Aspremont, and St{\'e}phane Mallat.
\newblock Phase recovery, {MaxCut} and complex semidefinite programming.
\newblock {\em Mathematical Programming}, 149:47--81, 2015.

\bibitem{Warnke2017UpperTails}
Lutz Warnke.
\newblock Upper tails for arithmetic progressions in random subsets.
\newblock {\em Israel Journal of Mathematics}, 221:317--365, 2017.

\bibitem{Zuo2003}
J.~M. Zuo, I.~Vartanyants, M.~Gao, R.~Zhang, and L.~A. Nagahara.
\newblock Atomic resolution imaging of a carbon nanotube from diffraction intensities.
\newblock {\em Science}, 300(5624):1419--1421, 2003.

\end{thebibliography}
\end{document}